\documentclass[9pt,amsfonts]{article}
\usepackage{graphicx,amssymb,eucal,mathrsfs}
\usepackage{latexsym,amsmath,amscd,epsfig,xcolor}
\usepackage[all]{xy}
\usepackage{ntheorem,tikz}
\usepackage{epsf,graphicx,pgf,etex,amscd,pstricks}

\usepackage[hypertexnames=false,backref=page,pdftex,
 	pdfpagemode=UseNone,
 	breaklinks=true,
 	extension=pdf,
 	colorlinks=true,
 	linkcolor=blue,
 	citecolor=red,
 	urlcolor=blue,
 ]{hyperref}

\newcommand{\Z}{\mathbb{Z}}

\newcommand{\C}{\mathbb{C}}

\newcommand{\del}{{\partial}}
\newcommand{\mcg}{{\rm Mod}}
\newcommand{\out}{{\rm Out}}

\newtheorem{main}{Theorem}

\newtheorem{prop}[main]{Proposition}
\newtheorem{lemma}[main]{Lemma}
\newtheorem{rem}[main]{Remark}
\newtheorem{cor}[main]{Corollary}

\theoremstyle{definition}

\newtheorem{defi}[main]{Definition}

\renewcommand{\thefootnote}{\alph{footnote}}

\title{Mapping class groups, multiple Kodaira fibrations, and ${\rm CAT}(0)$ spaces}
\date{November 2020}
\author{Claudio Llosa Isenrich$^{{\rm a}}$ and Pierre Py$^{{\rm b}}$}

\begin{document}

\maketitle

\begin{abstract} 
We study several geometric and group theoretical problems related to Kodaira fibrations, to more general families of Riemann surfaces, and to surface-by-surface groups. First we provide constraints on Kodaira fibrations that fiber in more than two distinct ways, addressing a question by Catanese and Salter about their existence. Then we show that if the fundamental group of a surface bundle over a surface is a ${\rm CAT}(0)$ group, the bundle must have injective monodromy (unless the monodromy has finite image). Finally, given a family of closed Riemann surfaces (of genus $\ge 2$) with injective monodromy $E\to B$ over a manifold $B$, we explain how to build a new family of Riemann surfaces with injective monodromy whose base is a finite cover of the total space $E$ and whose fibers have higher genus. We apply our construction to prove that the mapping class group of a once punctured surface virtually admits injective and irreducible morphisms into the mapping class group of a closed surface of higher genus.\end{abstract}

\tableofcontents

\footnotetext[1]{The author is grateful to the Max Planck Institute for Mathematics for its hospitality and financial support.}
\footnotetext[2]{Partially supported by the french project ANR AGIRA.}

\newpage

\renewcommand{\thefootnote}{\arabic{footnote}}

\section{Introduction}

A {\it Kodaira fibration} is a compact complex surface $X$ endowed with a holomorphic submersion onto a Riemann surface $\pi : X \to \Sigma$ which has connected fibers and is not isotrivial. The genus of the base and the fiber of $\pi$ must necessarily be greater than $1$. Since the construction of the first examples of such fibrations by Atiyah and Kodaira~\cite{Ati-69,Kod-67}, these complex surfaces have been widely studied, either from the point of view of complex geometry or from the point of view of the theory of surface bundles. See for instance~\cite{bregman,catanese,cr,caupol,chen,flapan1,flapan2,gdh,llr,Sal-15} for a few works on this topic. In this text we study several geometric and group theoretical problems related to these complex surfaces, to more general families of Riemann surfaces, and to surface-by-surface groups. 

Our work consists of three main parts. Firstly, we address the question whether Kodaira fibrations can fiber in more than two ways (Section \ref{sectionmf}). The question of the existence of such multiple fibrations has been raised independently by Catanese~\cite{catanese} and Salter~\cite{Sal-15}. We apply various tools from the theory of isolated singularities to provide strong constraints on the existence of such fibrations (Theorem~\ref{resomf}). In particular, we use a classical theorem of Mumford about normal singularities of complex surfaces~\cite{Mum-61}. Secondly, we discuss the properties of surface-by-surface groups which are ${\rm CAT}(0)$ (Section \ref{cat}). Using work of Monod on actions of direct products on ${\rm CAT}(0)$ spaces~\cite{monod}, we prove that a necessary condition for such a group to be ${\rm CAT}(0)$ is to have injective monodromy, provided the monodromy has infinite image (Theorem~\ref{theo-cat}). Thirdly, we introduce a technique that produces new families of Riemann surfaces with injective monodromy out of old ones (Sections \ref{htmcg} and \ref{geometricpart}). This relies on some old work of F.E.A. Johnson~\cite{j94}. One application of this is the construction of new examples of virtual injections between certain mapping class groups (Theorem~\ref{veofmcg}).

Before stating our results more precisely, we fix some notations and introduce additional definitions. Since both real and complex surfaces will appear throughout this work, we will use the following convention to distinguish between them.

\vspace{.4cm}

{\noindent \bf Convention.} If we do not explicitly state otherwise, then by a surface we always mean a real oriented surface. In particular, we will be careful to write complex surface when we refer to a space of complex dimension 2. 

As for the notations, if $F$ is a surface, we denote by $\mcg (F)$ its mapping class group and by $\mcg (F,\ast)$ the mapping class group of the surface $F$ with one marked point. Occasionally, we also write $S_g$ to denote the closed oriented surface of genus $g$. 

\vspace{.4cm}

Before discussing our first theorem, which concerns {\it multiple Kodaira fibrations}, we state the following notion of equivalence between fiber bundles. Although it is not completely standard, it is the natural one in the context of our work. 

\begin{defi}\label{def:distinctfibgeometric} 
Let $M$ be a smooth manifold and let $f_i:M\to B_i$, $i=1,2$, be two smooth submersions with connected fibers. We call them equivalent if the induced maps on fundamental groups $f_{i,\ast}:\pi_1(M)\to \pi_1(B_i)$ have the same kernel and we say that they are distinct otherwise.
\end{defi}

There are several classical results concerning the existence of distinct differentiable fibrations on a given $3$- or $4$-dimensional manifold. If $M$ is a hyperbolic $3$-manifold which fibers over the circle, then it is known that $M$ fibers in more than one way if and only if its first Betti number is greater than $1$, and if this is the case then $M$ fibers in infinitely many ways, by results of Thurston \cite{Thu-86}. In contrast, if $M$ is a smooth $4$-manifold it admits only finitely many distinct fibrations over a surface of genus greater than $1$ with fiber of genus greater than $1$. This is a result of Johnson~\cite{j99} (see also \cite{Sal-16}). Salter proved that the number of such distinct fibrations can be bounded above in terms of the Euler characteristic of $M$ and gave examples of $4$-dimensional manifolds where the number of distinct fibrations is arbitrarily large~\cite{Sal-15}. However, Salter's examples cannot admit a complex structure \cite{Sal-15}. Presently, there are examples of Kodaira fibrations which are known to fiber in only one way or which admit at least two distinct fibrations. For some examples,we know that the number of fibrations is exactly two, see for instance the article~\cite{chen} which proves that some of the classical Atiyah-Kodaira examples fiber in {\it exactly} two distinct ways.  But there is currently no example of a Kodaira fibration which is known to fiber in more than two distinct ways. This led Catanese \cite{catanese} and Salter \cite{Sal-15} to independently ask the question whether for a Kodaira fibration the number of distinct fibrations can be greater than $2$. 

We prove here the following result which can be seen as a restriction on potential examples of complex surfaces admitting three distinct Kodaira fibrations.

\begin{main}\label{resomf} Let $X$ be a compact complex surface. Suppose that $X$ admits three distinct Kodaira fibrations $p_i : X\to \Sigma_{i}$, $1\le i \le 3$. Then the image of the induced morphism
$$\pi_{1}(X)\to \pi_{1}(\Sigma_{1})\times \pi_{1}(\Sigma_{2})\times \pi_{1}(\Sigma_{3})$$
is of finite index in $\pi_{1}(\Sigma_{1})\times \pi_{1}(\Sigma_{2})\times \pi_{1}(\Sigma_{3})$. 
\end{main}

The reader will find in Remark~\ref{uneseulesubmersion} a slight reinforcement of this theorem. Note that Kodaira fibrations are always projective, hence K\"ahler~\cite[Ch. V]{ccs}. If $X$ is a compact K\"ahler manifold, the number of distinct surjective holomorphic maps (possibly with critical points) with connected fibers from $X$ to closed hyperbolic Riemann surfaces (or orbifolds) is always finite, see for instance~\cite{CorSim-08,Del-08}. This raises the question of finding upper bounds on the number of such maps for a fixed choice of a compact K\"ahler manifold. Catanese and Salter's question fits into that context. Holomorphic maps from compact K\"ahler manifolds to closed Riemann surfaces play a central role in the study of K\"ahler groups, see~\cite{ABCKT-95}. Recently there has been considerable progress in understanding holomorphic maps from compact K\"ahler manifolds to direct products of closed Riemann surfaces (see e.g. \cite{DelGro-05,DelPy-16,Llo-17,Llo-18-I}). Our proof of Theorem~\ref{resomf} will rely on these ideas and notably on the work of the first author~\cite{Llo-18-I}, as well as on some tools from the theory of isolated singularities.

Theorem~\ref{resomf} raises a natural question. Assume that $X$ is a compact complex surface admitting three distinct Kodaira fibrations $p_i : X\to \Sigma_i$ ($1\le i \le 3$). Is it true that the image of $X$ in the product $\Sigma_1\times \Sigma_2 \times \Sigma_3$ is an ample divisor?  


 Before turning to the presentation of the second and third part of our work, we state some group theoretical definitions. From now on, by a {\it surface group} we will mean the fundamental group of a closed oriented surface of genus greater than $1$. A {\it surface-by-surface group} is a group $G$ which admits a normal subgroup $R\lhd G$ such that both $R$ and $G/R$ are isomorphic to surface groups. Recall also that the group  $\out(G)$ of {\it outer automorphisms} of a group $G$ is defined as the quotient $\out(G)={\rm Aut}(G)/{\rm Inn}(G)$ of the group of automorphisms of $G$ by the group of inner automorphisms. If $R \lhd G$ is a normal subgroup of an arbitrary group $G$, there is a natural morphism 
\begin{equation*}
\varrho : G/R\to \out(R)
\end{equation*}
called the {\it monodromy morphism} of the extension:
\begin{equation}\label{extensiondebase}
\xymatrix{0\ar[r] & R \ar[r] & G \ar[r] & G/R \ar[r] & 0.}
\end{equation}
The properties of the monodromy of extensions as above with $R$ a surface group, in particular its injectivity, will be our main object of focus in the rest of this introduction.  

We now turn to our second main result. We observe that many examples of Kodaira fibrations admit ${\rm CAT}(0)$ metrics (see~\cite{bh} for basic facts about ${\rm CAT}(0)$ spaces). This is the case for all double etale Kodaira fibrations, in the terminology of~\cite{cr}. As we will explain in Section~\ref{cat}, this follows from standard results on ramified coverings and nonpositive curvature but does not seem to be widely known (see however~\cite{stadler}). It is then natural to ask which Kodaira fibrations have a fundamental group which is a ${\rm CAT}(0)$ group. We will see that a necessary condition for this is that the monodromy of the fibration is injective, provided that it has infinite image. This result applies in a purely group theoretical context, namely: 

\begin{main}\label{theo-cat} Let $G$ be a surface-by-surface group with infinite monodromy. We fix a normal subgroup $R\lhd G$ isomorphic to a surface group with $G/R$ isomorphic to a surface group. If the group $G$ is ${\rm CAT}(0)$ then the monodromy $G/R\to \out(R)$ is injective. 
\end{main}

The proof of the above theorem relies on the fact that if the monodromy is not injective, then $G$ must contain a nontrivial direct product. We then use a splitting theorem of Monod~\cite{monod} for actions of product groups on spaces of nonpositive curvature. It would be interesting to obtain further restrictions on surface-by-surface groups which are ${\rm CAT}(0)$ (or, narrowing our focus, on fundamental groups of Kodaira fibrations which are also ${\rm CAT}(0)$ groups).


We now turn to the last part of our work. Mirroring Definition~\ref{def:distinctfibgeometric}, we will consider two extensions of the same group $G$ as in~\eqref{extensiondebase} as distinct if the corresponding normal subgroups are distinct. In~\cite{j94}, F.E.A. Johnson proved the following theorem.

\begin{main}[Johnson]\label{theoremedejohnson} Let $G$ be a group which can be realized in two distinct ways as a surface-by-surface group. Assume that at least one of these realizations has infinite monodromy. Then the monodromy homormophism associated to any realization of $G$ as a surface-by-surface group is injective.
\end{main}

In what follows, we apply Johnson's theorem ``in family" to obtain new examples of groups containing a surface group as a normal subgroup and for which the monodromy is injective. To state our first result in this direction, we need to consider the higher dimensional {\it iterated Kodaira fibrations}, which we now define. The definition is made by induction, with the convention that a $2$-dimensional iterated Kodaira fibration is simply a compact complex surface which admits the structure of a Kodaira fibration.

\begin{defi}\label{itKodfib}
We say that an $n$-dimensional compact complex manifold $X$ is an iterated Kodaira fibration if there exists a holomorphic submersion $\pi : X \to Y$ with connected fibers which is not isotrivial, where $Y$ is an $(n-1)$-dimensional iterated Kodaira fibration.
\end{defi}

In the following, the expression {\it Kodaira fibration} will be reserved for complex surfaces. When talking about higher dimensional iterated Kodaira fibrations, we will always mention their dimension. If $\pi : X \to Y$ is as in the definition, one has a natural representation $\pi_{1}(Y)\to \mcg(F)$ where $F$ is a fiber of $\pi$. This is the {\it monodromy representation} of the fibration. It can be (equivalently) defined either in geometric terms or by considering the exact sequence of fundamental groups induced by $\pi$ and applying the algebraic definition given above. It is classical that one can construct iterated Kodaira fibrations in all dimensions, see for instance Miller's article~\cite{miller1984}. We prove here: 
 
\begin{main}\label{injm} For each $n\ge 2$, there exists an $n$-dimensional iterated Kodaira fibration $\pi : X \to Y$ with fiber $F$, such that the monodromy representation $\pi_{1}(Y)\to \mcg (F)$ is injective. 
\end{main}

A key ingredient in our proof of Theorem~\ref{injm} is Proposition~\ref{propositioninjg} which is a fibered version of Johnson's Theorem~\ref{theoremedejohnson}. Our arguments actually give more precise results. For instance, when $n=3$ in the previous theorem, we obtain the following (see Section \ref{candfi} for further results). 

\begin{main}\label{injm2}
If $X$ is a Kodaira fibration with injective monodromy, then there exists a family of closed Riemann surfaces $Z\to X'$ above some finite covering space $X'$ of $X$ whose monodromy is injective.
\end{main}

We do not know examples of Kodaira fibrations with injective monodromy which fiber in only one way. Hence the only current applications of Theorem \ref{injm2} are to double Kodaira fibrations, i.e. Kodaira fibrations which fiber in two distinct ways. However, we have chosen to formulate the theorem in the above way, since it might apply more generally. Observe also that, precisely because of Johnson's theorem mentioned above, one can talk of the property that a Kodaira fibration has injective monodromy independently of the choice of the fibration. Indeed, if the fibration is not unique, all fibrations must have injective monodromy by Johnson's result.

By generalizing Miller's construction~\cite{miller1984}, we also show that starting from any holomorphic submersion $\pi: Z\to B$ of complex manifolds defining a family of closed Riemann surfaces (of genus greater than $1$) with injective monodromy, one can construct new families of Riemann surfaces with injective monodromy whose base is a finite covering space of $Z$. We refer to Theorem~\ref{basearbitraire} for a more detailed statement. Here we only state the following consequence of that theorem.

\begin{main}\label{veofmcg} Let $g\ge 2$. Then there exists a finite index subgroup $\Gamma < \mcg (S_g,\ast)$ and an injective morphism $\phi : \Gamma \to \mcg (S_{2+8(g-1)})$ with the property that the action of $\phi(\Gamma)$ on the space of isotopy classes of simple closed curves on $S_{2+8(g-1)}$ does not have any finite orbit. 
\end{main}

Recall that a subgroup of the mapping class group of a closed surface is {\it irreducible} if it does not fix any isotopy class of simple closed curve. Our result then says that any finite index subgroup of $\phi (\Gamma)$ is irreducible. We also mention that our morphism comes with an equivariant holomorphic map between the corresponding Teichm\"uller spaces, see Remark~\ref{appholoeq}. 

To motivate Theorem~\ref{veofmcg}, we recall that various authors have studied morphisms between different mapping class groups, with a particular focus on their rigidity properties. The general idea, partly inspired by Margulis superrigidity and by the comparison of mapping class groups with lattices, is that nontrivial morphisms between mapping class groups should be induced by geometric constructions obtained by ``manipulation of surfaces". We refer the reader to~\cite{arso,arleso,arso-rigidity,bridson} for more details. The morphism constructed in Theorem \ref{veofmcg} seems to be obtained from a more complicated manipulation of surfaces than the earlier examples. It was previously known that the group $\mcg (S_g,\ast)$ embeds into the group $\mcg (S_{2g})$, see~\cite[\S 3.3]{arso-rigidity}. However, the morphism constructed in~\cite{arso-rigidity} is reducible. In contrast, the morphism we build is irreducible but is only defined on a finite index subgroup and its range is a mapping class group of bigger genus. We will go back to this topic and prove Theorem~\ref{veofmcg} in Section~\ref{sectionmcg}.

\vspace{.4cm}

The text is organized as follows. Section~\ref{sectionmf} studies complex surfaces admitting three distinct Kodaira fibrations; it contains the proof of Theorem~\ref{resomf}. Section~\ref{cat} contains the proof of Theorem~\ref{theo-cat}. Section~\ref{htmcg} contains some group theoretical results and notably Proposition~\ref{propositioninjg} which gives a way to build extensions by surface groups with injective monodromy. This material is needed in Section~\ref{geometricpart} which contains some geometric constructions of families of Riemann surfaces and includes the proofs of Theorems~\ref{injm}, \ref{injm2} and \ref{veofmcg}. Finally, Section~\ref{furtherremarks} contains a few more observations about the fundamental groups of Kodaira fibrations, including an answer to a question of Bregman about their residual nilpotence. 

The three main parts of this work, namely Section~\ref{resomf}, Section~\ref{cat} and Sections~\ref{htmcg} and~\ref{geometricpart} together are essentially independent from each other and can be read independently. 

\vspace{.4cm}

{\bf Acknowledgements.} We would like to thank Patrick Popescu-Pampu for a useful correspondence, Domingo Toledo for his comments on a preliminary version of this article, as well as the mathematical institute of the UNAM in Mexico city, where our collaboration started. We also would like to thank Nick Salter for pointing out Morita's work~\cite{morita1,morita2} to us as well as the referee for their constructive comments on the text.


\section{Multiple fibrations}\label{sectionmf}

This section is devoted to the proof of Theorem~\ref{resomf}. In Section~\ref{secLemmas}, we collect some simple lemmas and observations about complex surfaces admitting several distinct Kodaira fibrations, as well as results due to Bridson, Howie, Miller and Short on the one hand, and the first author of the present article on the other hand. Section~\ref{isosing} collects a few results about isolated singularities of complex surfaces and of holomorphic maps. Finally, Section~\ref{potmtttt} contains the proof of Theorem~\ref{resomf}.

\subsection{Preliminary results}\label{secLemmas}

As we shall see at the very end of Section~\ref{potmtttt}, besides Theorem~\ref{resomf}, we can also obtain restrictions on compact complex surfaces admitting several fibrations onto Riemann surfaces even if we allow fibrations which are not submersions, as soon as {\it at least one of the fibrations is a submersion}. For this reason, we introduce fibrations onto Riemann surfaces in general and state our preliminary lemmas in this more general context.

\begin{defi}\label{def-fibration} A fibration from a compact complex manifold $X$ onto a Riemann surface $\Sigma$ is a surjective holomorphic map with connected fibers. 
\end{defi}

If $f : X \to \Sigma$ is a fibration and if $F$ is a singular fiber of $f$, one can define its multiplicity $m(F)$, see e.g.~\cite{catanesekeumoguiso,Del-08} for the definition. The Riemann surface $\Sigma$ then inherits a natural orbifold structure by assigning to the image $p=f(F)$ of each critical fiber $F$ the multiplicity $m(F)$. The corresponding orbifold fundamental group is denoted by $\pi_{1}^{orb}(\Sigma)$. The natural morphism from the fundamental group of $X$ to that of $\Sigma$ lifts to a morphism $\pi_{1}(X)\to \pi_{1}^{orb}(\Sigma)$ that we still denote by $f_{\ast}$ (where it is implicit that the orbifold structure is induced by $f$). The kernel of this morphism is generated by the image of a generic fiber of $f$, hence is finitely generated. For all these facts we refer the reader to~\cite{catanesekeumoguiso,Del-08} for instance. Similarly, we follow~\cite{Del-08} to state the next definition.

\begin{defi}\label{distinctfib} Let $X$ be a compact complex manifold. Two fibrations $f_{i} :X \to \Sigma_{i}$ ($i=1, 2$) onto some Riemann surfaces are equivalent if they have the same fibers, the same singular fibers and if the multiplicities of their singular fibers are the same. We say that $f_{1}$ and $f_{2}$ are distinct otherwise. 
\end{defi}

Of course if $f_1$ and $f_2$ are submersions, $f_1$ and $f_2$ are equivalent if and only if they have the same fibers; there is no need to consider singular fibers. Observe that two such submersions are equivalent in the sense of Definition \ref{distinctfib} if and only if they are equivalent in the sense of Definition \ref{def:distinctfibgeometric}. To avoid cumbersome notations, from now on and until the end of Section~\ref{potmtttt}, we will assume that all fibrations that we encounter have no multiple fibers. This includes maps which can have critical points and this allows us to deal exclusively with surface groups (instead of orbifold surface groups). Since one can always reduce to this situation by taking finite coverings, this is not a big restriction. We will only go back to the most general case (including possibly multiple fibers) in Remark~\ref{uneseulesubmersion}.

We also point out that in this section, the word {\it fibration} will refer to Definition~\ref{def-fibration}. We will use the expression {\it Kodaira fibration} to designate a fibration defined by a submersion. We now turn to some preliminary results. We start with the following elementary lemma, whose proof will be omitted. 

\begin{lemma}\label{lem2D}
 Let $X$ be a compact complex manifold which admits $r\geq 2$ distinct fibrations $f_i: X \to \Sigma_i$, $1\leq i \leq r$. For each point $s_{i}\in \Sigma_{i}$, denote by $F_{i,s_i}$ the fiber $f_{i}^{-1}(s_{i})$. Then for $i\neq j$ and $s_j\in \Sigma_j$ the restriction $f_j|_{F_{j,s_j}}: F_{j,s_j} \to \Sigma_i$ is surjective. In particular, the holomorphic map $\left(f_i,f_j\right):X\to \Sigma_i\times \Sigma_j$ is surjective for $1\leq i<j\leq r$.\end{lemma}
 
 \begin{cor}\label{etablirvsp}
  The image of the map $\pi_{1}(X)\to \pi_{1}(\Sigma_{i})\times \pi_{1}(\Sigma_{j})$ induced by $(f_{i},f_{j})$ has finite index in $\pi_{1}(\Sigma_{i})\times \pi_{1}(\Sigma_{j})$. 
\end{cor}
{\it Proof.} This is a general property of surjective holomorphic maps between compact complex manifolds.\hfill $\Box$

\vspace{.4cm}

\noindent To prove Theorem~\ref{resomf}, we will need to study the image of the map $$X \to \Sigma_{1}\times \Sigma_{2}\times \Sigma_{3}$$ 
induced by three distinct Kodaira fibrations on the same complex surface. The following lemma shows that, if smooth, such an image is well-understood.  

\begin{lemma}\label{lemSmooth}
 Let $X$ be a complex surface and let $f_i: X\to \Sigma_{i}$, $1\leq i \leq 3$, be distinct fibrations. Denote by $f=\left(f_1,f_2,f_3\right):X\to \Sigma_1\times \Sigma_2 \times \Sigma_3$ the product map. If the image $Y=f(X)$ is smooth then it admits itself three distinct fibrations induced by the projections to the $\Sigma_i$. Furthermore, if $f_{i_{0}}$ was a submersion, the induced map $Y\to \Sigma_{i_{0}}$ is also a submersion. 
\end{lemma}

\noindent {\it Proof.} We assume that $Y$ is smooth and denote by $p_i: \Sigma_1 \times \Sigma_2 \times \Sigma_3\to \Sigma_{i}$ the projection onto the factor number $i$ and by $g_i=p_i|_Y$ its restriction to $Y$. Then $f_i$ decomposes as
 \[
  \xymatrix{ X \ar[r]^f \ar[rd]_{f_i} & Y\ar[d]^{g_i}\\ & \Sigma_i.}
 \]
Since $X$, $Y$ and $\Sigma_i$ are smooth, the chain rule implies that the map $g_i$ is a submersion if $f_i$ was so. The fibers of $Y\to \Sigma_i$ are connected, since they are images of the fibers of $f_i$. This implies that $Y$ admits three distinct fibrations.\hfill $\Box$

We now observe that distinct fibrations on a compact complex surface $X$ induce distinct fibrations on any finite covering space of $X$. More precisely, let $f_i:X\to \Sigma_i$, $1\leq i \leq r$, be distinct fibrations, and let $h: X_0\to X$ be a finite covering. Denote by $p_i: \Sigma'_i\to \Sigma_i$ the coverings corresponding to the subgroups $f_{i,\ast}(\pi_1(X_0))\leq \pi_1(\Sigma_i)$. Then there are induced holomorphic maps $f_{i,0}:X_0\to \Sigma'_i$ making the diagram
\begin{equation}
\label{eq:commDiag1}
\begin{gathered}
\xymatrix{ X_0\ar[r]^{h}\ar[d]_{f_{i,0}}& X\ar[d]^{f_i} \\ \Sigma'_i \ar[r]^{p_i} & \Sigma_i} 
\end{gathered}
\end{equation}
commutative. We then have:

\begin{lemma}\label{inducedffgg}
The maps $f_{i,0}: X_0\to \Sigma'_i$ ($1\le i \le r$) are pairwise distinct fibrations. If $f_i$ was a submersion then so is $f_{i,0}$. 
 \label{lemFinCover}
\end{lemma}
{\it Proof.} It is clear that the map $f_{i,0}$ has connected fibers and is a submersion if $f_i$ was one. The fact that the $r$ fibrations $f_{i,0}: X_0\to \Sigma'_i$ are distinct is an immediate consequence of Lemma \ref{lem2D} applied to the $f_i$ and \eqref{eq:commDiag1}. Indeed, they imply that for $j\neq i$ the restriction of $f_{i,0}$ to any fiber of $f_{j,0}$ surjects onto $\Sigma'_i$.\hfill $\Box$

\vspace{.4cm}

We now introduce some group theoretical notions which will be needed to formulate the results from~\cite{BriHowMilSho-13} that we will use. If $G_{1}\times \cdots \times G_{n}$ is a direct product of arbitrary groups $(G_{i})_{1\le i \le n}$, and if $1\le i<j\le n$ are distinct indices, we denote by $p_{ij} : G_{1}\times \cdots \times G_{n}\to G_{i}\times G_{j}$ and $p_{i} : G_{1}\times \cdots G_{n}\to G_{i}$ the natural projections. 

\begin{defi} A subgroup $H < G_{1} \times \cdots \times G_{n}$ of a direct product is said to {\it virtually surject onto pairs} if for any pair of indices $i<j$ the group $p_{ij}(H)$ has finite index in $G_{i}\times G_{j}$. 
\end{defi}

\begin{defi} A subgroup $H < G_{1} \times \cdots \times G_{n}$ of a direct product is a {\it subdirect product} if $p_{i}(H)=G_{i}$ for every $i\in \{1, \ldots , n\}$. It is full if $H \cap G_{i} \neq\{1\}$ for every index $i\in \{1, \ldots , n\}$ (where $G_{i}$ is naturally embedded in the direct product).
\end{defi}

We can now state the following theorem (see Theorem A in~\cite{BriHowMilSho-13}):

\begin{main}[Bridson, Howie, Miller, Short] Let $G_{1}, \ldots , G_{n}$ be finitely presented groups. Let  $H < G_{1} \times \cdots \times G_{n}$ be a subgroup which virtually surjects onto pairs. Then $H$ is finitely presented. 
\end{main}

This is part of a long series of results of these authors about subgroups of direct products in general, and more specifically subgroups of direct products of limit groups. We refer the reader to~\cite{bhms2002,bhms2009} for more results in this vein. Using this and the previous observation, we shall prove:

\begin{prop}\label{propCoab}
 Assume that $f=(f_1,f_2,f_3): X \to \Sigma_1\times \Sigma_2 \times \Sigma_3$ is a compact complex surface admitting three distinct fibrations $f_i$. Then $f_{\ast}(\pi_{1}(X))< \pi_1(\Sigma_1)\times \pi_{1}(\Sigma_2)\times \pi_1(\Sigma_{3})$ is a finitely presented full subdirect product. 
 \end{prop} 
 {\it Proof.} The fact that $f_{\ast}(\pi_{1}(X))$ is finitely presented follows from Corollary~\ref{etablirvsp}, combined with Bridson, Howie, Miller and Short's theorem mentioned above. The fact that it is subdirect is clear and we must justify that it is full. Although this follows from the results in~\cite[\S 5]{Llo-18-I} (see Lemma 5.3 there), we give the argument here as it is quite simple. Assume by contradiction that $f_{\ast}(\pi_{1}(X))$ is not full. Up to reordering the factors, we assume that $f_{\ast}(\pi_{1}(X)) \cap \pi_{1}(\Sigma_{1})=\{1\}$. This means that the projection $\pi_1(\Sigma_1)\times \pi_1(\Sigma_2)\times \pi_{1}(\Sigma_{3})\to \pi_1(\Sigma_2)\times \pi_{1}(\Sigma_{3})$ induces an isomorphism from $f_{\ast}(\pi_{1}(X))$ onto a finite index subgroup of $\pi_1(\Sigma_2)\times \pi_{1}(\Sigma_{3})$, which is nothing else than $(f_2,f_3)_{\ast}(\pi_{1}(X))$. This implies that the morphism 
 $$f_{1,\ast} : \pi_{1}(X) \to \pi_{1}(\Sigma_{1})$$
 factors through the map $(f_2,f_3)_{\ast}$, i.e. there exists a morphism $$\varphi : (f_2,f_3)_{\ast}(\pi_{1}(X))\to \pi_{1}(\Sigma_1)$$ such that 
 \begin{equation}\label{factorupl}
 f_{1,\ast}=\varphi\circ (f_2,f_3)_{\ast}.
 \end{equation}
 By Corollary~\ref{etablirvsp}, the group $(f_{2},f_{3})_{\ast}(\pi_{1}(X))$ has finite index in $\pi_{1}(\Sigma_{2})\times \pi_{1}(\Sigma_{3})$, hence contains a finite index subgroup of the form $\pi_{1}(\Sigma'_2)\times \pi_{1}(\Sigma'_{3})$ for some finite covering spaces $\Sigma_{i}'\to \Sigma_{i}$. Let $X_{0}\to X$ be the finite covering space such that $$\pi_{1}(X_{0})=(f_{2},f_{3})_{\ast}^{-1}(\pi_{1}(\Sigma'_2)\times \pi_{1}(\Sigma'_{3})).$$  Equation~\eqref{factorupl} then implies that $f_{1,\ast}(\pi_{1}(X_{0}))=\varphi (\pi_1(\Sigma'_2)\times \pi_1(\Sigma'_3))$. Since this group has finite index in $\pi_{1}(\Sigma_{1})$ and since a surface group cannot be generated by two nontrivial commuting subgroups we must have $\varphi (\pi_{1}(\Sigma'_2))=\{1\}$ or $\varphi (\pi_{1}(\Sigma'_3))=\{1\}$. Let us assume that $\varphi (\pi_1(\Sigma'_2))=\{1\}$. This means that the restriction of $f_{1,\ast}$ to $\pi_{1}(X_{0})$ factors through the map $f_{3,\ast}$. But this contradicts the fact that $f_1$ and $f_3$ induce distinct fibrations on $X_{0}$ (see Corollary~\ref{etablirvsp} and Lemma~\ref{inducedffgg}).\hfill $\Box$
 
 \vspace{.4cm}

Before stating the last result of this section we recall some developements from the last fifteen years, around the construction of new examples of K\"ahler groups. In~\cite{dps}, Dimca, Papadima and Suciu built new examples of K\"ahler groups using the following construction. They consider a finite number of Riemann surfaces $\Sigma_1, \ldots , \Sigma_n$ each admitting a ramified covering $q_i : \Sigma_i \to E$ of degree two over the same elliptic curve $E$. They then prove that if $n\ge 3$, the fundamental group of a smooth fiber of the map $q_1 +\cdots + q_n$ from the direct product of the $\Sigma_{i}$'s to $E$ injects into the product $\pi_{1}(\Sigma_1)\times \cdots \times \pi_{1}(\Sigma_{n})$. Moreover the corresponding group has some exotic finiteness properties. This construction was pushed further by the first author~\cite{Llo-fourier,Llo-17} to build more examples of K\"ahler groups. In a related direction, the article~\cite{Llo-18-I} studies images of K\"ahler groups in a direct product of surface groups when the morphism is induced by a holomorphic map. In particular, it studies the following situation. Consider a compact K\"ahler manifold $X$ and let $p_i : X \to \Sigma_i$ ($1\le i \le 3$) be surjective holomorphic maps with connected fibers onto  some Riemann surfaces. Let $p=(p_1,p_2,p_3) : X \to \Sigma_1 \times \Sigma_2 \times \Sigma_3$. What can one say about the subgroup $p_{\ast}(\pi_{1}(X)))$? The following result from~\cite{Llo-18-I} answers this question (note that~\cite{Llo-18-I} also contains further results around the same question). 

\begin{main}\label{clfsp} Assume that the image $p_{\ast}(\pi_{1}(X))$ of $\pi_{1}(X)$ in $\pi_{1}(\Sigma_1)\times \pi_{1}(\Sigma_2)\times \pi_1(\Sigma_3)$ is a finitely presented full subdirect product of {\it infinite index}. Then there exist finite covering spaces $\Sigma'_i \to \Sigma_i$ ($1\le i \le 3$), an elliptic curve $E$ and some ramified coverings $q_i : \Sigma'_i \to E$ such that the kernel of the map $$(q_1+q_2+q_3)_{\ast} : \pi_1(\Sigma'_1)\times \pi_1(\Sigma'_2)\times \pi_1(\Sigma'_3)\to \pi_1(E)$$ is contained in $p_{\ast}(\pi_{1}(X))$ as a subgroup of finite index. 

Moreover, let $X_0\to X$ be the finite covering space corresponding to the subgroup $p_{\ast}^{-1}(\pi_1(\Sigma'_1)\times \pi_1(\Sigma'_2)\times \pi_1(\Sigma'_3))$. Then the image of $X_0$ in $\Sigma'_1 \times \Sigma'_2 \times \Sigma'_3$ coincides with a (possibly singular) fiber of the map $q_1+q_2+q_3$.    
\end{main}

Note that this theorem is not stated in this exact form in~\cite{Llo-18-I} but it follows from Theorems 1.1 and 3.1 there.

\subsection{Isolated singularities}\label{isosing}

In this section we recall some classical facts about isolated singularities of complex spaces and of holomorphic maps. Our introduction will be based on \cite{Loo-84} and we refer the reader to that book for a detailed introduction to the subject. For an introduction to complex analytic spaces, we refer the reader to~\cite{gr}. 

Let $Y$ be a complex analytic space and let $y_{0}\in Y$. We shall always assume that $Y-\{y_{0}\}$ is smooth. If $Y$ is not smooth at $y_{0}$, we say that $y_{0}$ is an {\it isolated singular point} of $Y$. Note that most of the results from~\cite{Loo-84} that we will use are valid in greater generality, but this simple case is enough for our purpose. If $U$ is a neighborhood of $y_0$ in $Y$, a function $r: U\to \left[0,\infty\right)$ \textit{defines} the point $y_0$ if $r^{-1}(0)= \left\{y_{0}\right\}$ and $r$ is real analytic (this means that $U$ can be realized as an analytic set in some open set $W\subset \C^N$ and that $r$ extends to a real analytic function on a neighborhood of $U$ inside $W$). Note that such maps $r$ always exists for a suitable choice of $U$. \vspace{.4cm}

We now fix such a function $r$. Up to shrinking $Y$, we can and do assume that $r$ is defined on $Y$. Since $r$ is analytic, there is $\varepsilon_0>0$ such that $0$ is the only critical value of $r$ in $\left[0,\varepsilon_0\right]$. We consider a holomorphic map $f : Y \to \C$ such that $f(y_{0})=0$. We assume that $f$ is a submersion on $Y-\{y_{0}\}$. For $\varepsilon>0$ and for a subset $M\subset \C$, we will use the notation $Y_{M, r\leq \varepsilon}$, $Y_{M, r=\varepsilon}$, $Y_{M, r<\varepsilon}$ to denote respectively the sets
\[
 \left\{y\in Y\mid f(y)\in M,~~ r(y)\leq \varepsilon\right\},\mbox{~~~} \left\{y\in Y\mid f(y)\in M,~~ r(y) = \varepsilon\right\},
\] 
  \[\mbox{~~~} \left\{y\in Y\mid f(y)\in M,~~ r(y)< \varepsilon\right\}.\]

\noindent In what follows $D$ will be a closed disc centered at the origin in $\C$ and we will focus on sets of the form $Y_{D,r\le \varepsilon}$ which are schematically pictured in Figure 1.

\begin{center}
\begin{tikzpicture}
\draw [thick] (-.295,2) -- (.25,1.026) ;
\draw [thick] (-.304,1.8) -- (.128,1.028) ;
\draw [thick] (-.325,1.6) -- (0,1.019) ;
\draw [thick] (-.358,1.428) -- (-.138,1.03486) ;
\draw [thick] (-.407,1.26) -- (-.28611,1.0441) ;
\draw [thick] (-.31,2.27) -- (.375,1.045905) ;
\draw [thick] (-.347,2.52) -- (.483,1.05466) ;
\draw [thick] (-.345,2.698) -- (.565,1.07183) ;
\draw [thick] (-.212,2.712) -- (.688,1.1037) ;
\draw [thick] (-.05,2.714) -- (.635,1.489905) ;
\draw [thick] (.145,2.718) -- (.618,1.872749) ;
\draw [thick] (.315,2.702) -- (.635,2.13016) ;
\draw [thick] (.43,2.685) -- (.657,2.279351) ;
\draw [thick] (.61,2.65) -- (.7,2.48917) ;
\draw (0,1.87) ellipse (1.55cm and .85cm) ;
\draw (-1,1.8) node {{\small $Y_{r\le \varepsilon}$}} ;
\draw [->] (0,.9) -- (0,.5) ;
\draw (.45,.75) node {{\small $f$}} ;
\draw (0,0) circle (.25) ;
\draw (.96,3.1) arc (150:205:2.55) ;
\draw (-0.62,0.8) arc (-30:30:2.4) ;
\draw (0.1,3) node {{\small $Y_{D}$}} ;
\draw (0,0) node {{\small $\cdot$}} ;
\draw (0,-.7) node {Figure 1: the set $Y_{D,r\le \varepsilon}$} ;
\end{tikzpicture}
\end{center}

\begin{defi}
 For $\varepsilon\in \left(0,\varepsilon_0\right]$ we will say that the restriction $f: Y_{D,r\leq \varepsilon} \to D$ is a \textit{good proper representative} of $f$ if $f|_{r=\varepsilon}$ is a submersion at all points of $Y_{D,r=\varepsilon}$.
\end{defi}

One can always find arbitrarily small $\varepsilon$ and $D$ such that $f: Y_{D,r\leq \varepsilon} \to D$ is a good proper representative. The significance of this notion stems from the following general result.

\begin{main}[Theorem 2.8 in \cite{Loo-84}]\label{thmIsolsing}
 For a good proper representative $f: Y_{D,r\leq \varepsilon}\to D$ the following holds:
 \begin{enumerate}
\item $f$ is proper and $f: Y_{D,r=\varepsilon}\to D$ is a trivial smooth fiber bundle;
\item $f: \left( Y_{D\setminus \{0\},r\leq \varepsilon}, Y_{D\setminus \{0\}, r=\varepsilon}\right) \to D\setminus \{0\}$ is a smooth fiber bundle pair.
\end{enumerate}
\end{main}

\begin{defi}
 For $s\neq 0$, we call $Y_{s,r\leq \varepsilon}$ the Milnor fiber of $f|_{ Y_{D,r\leq \varepsilon}}$. 
\end{defi}

A particularly important result for our purposes is:
\begin{lemma}[Lemma 2.10 in~\cite{Loo-84}]
\label{lemCone}
Let $f: Y_{D,r\leq \varepsilon} \to D$ be a good proper representative. There is $\eta_0>0$ such that for any disc $D_{\eta}\subset D$ of radius $\eta_0\geq \eta >0$ the variety $Y_{D_{\eta},r\leq \varepsilon}$ is homeomorphic to the cone over $\del(Y_{D_{\eta},r\leq \varepsilon})$.
\end{lemma}

\begin{rem}\label{remarklink} Although this is not stated in~\cite{Loo-84}, we observe that $\del(Y_{D_{\eta},r\leq \varepsilon})$ is homeomorphic to the link of the singular point $y_{0}\in Y$. This follows readily from the arguments presented there (see~\cite[\S 2.A]{Loo-84} for the definition of the link of an isolated singular point). 
\end{rem}

In what follows, we always assume that whenever we choose a good proper representative as before, then the radius of $D$ is small enough such that $Y_{D,r\leq \varepsilon}$ is homeomorphic to the cone over its boundary.

\vspace{.4cm}

In the next section, we will need to rule out the existence of isolated singular points for certain (possibly singular) complex surfaces embedded in a product of three Riemann surfaces. The key tool for this will be the following proposition which considers isolated singular points on a normal complex analytic surface.   

\begin{prop}\label{lemMilnDiscSmooth}
Assume that $Y$ is normal of dimension $2$. Let $f:Y_{D,r\leq \varepsilon}\to D$ be a good proper representative. Assume that the Milnor fiber of $f$ is a disc. Then $Y$ is smooth at $y_{0}\in Y$.
\end{prop}
{\it Proof.} Since the Milnor fiber of $f$ is a disc, the map $Y_{\del(D),r\leq \varepsilon}\to \del (D)$ is a disc bundle over the circle. Topologically the only orientable disc bundle over the circle is the trivial bundle, implying that $Y_{\del (D), r\leq \varepsilon}$ is a solid torus with meridian $\mu=Y_{s,r=\varepsilon}$ for $s\in \del(D)$. Note that, by definition of the Milnor fiber, $\mu$ coincides with the longitude of the solid torus $Y_{D,r=\varepsilon}$, since it is a fiber of the locally trivial fibration $Y_{D,r=\varepsilon}\to D$. In particular, $\del(Y_{D,r\leq \varepsilon})$ is obtained by gluing two solid tori identifying the meridian of the first torus with the longitude of the second torus. Since this uniquely determines the resulting 3-manifold, we deduce that $\del(Y_{D,r\leq \varepsilon})$ is a topological 3-sphere. By Lemma \ref{lemCone}, $Y_{D,r\leq \varepsilon}$ is homeomorphic to a cone over this 3-sphere. This implies that there is a neighbourhood of $y_0$ in $Y$ which is a topological 4-manifold. It then follows from a famous theorem of Mumford that $y_0$ is actually a smooth point of $Y$~\cite[p.1]{Mum-61}.\hfill $\Box$

\subsection{Multiple fibrations and singularities}\label{potmtttt}

Consider a Kodaira fibration $X$ that admits three distinct fibrations $f_i : X \to \Sigma_i$ ($1\le i \le 3$) in the sense of Definition~\ref{def-fibration}. We let $f$ be the product map:
$$f=(f_1,f_2,f_3): X \to \Sigma_1\times \Sigma_2 \times \Sigma_3.$$
\noindent We assume that $f_{1}$ is a submersion. Let $Y:=f(X)$. It is a purely $2$-dimensional irreducible analytic subspace of $\Sigma_{1}\times \Sigma_{2}\times \Sigma_{3}$. We will prove the following proposition.

\begin{prop}\label{normalisation} The normalization $\widetilde{Y}$ of $Y$ is smooth. In particular $Y$ does not have any isolated singular point.  
\end{prop}

Let us explain why this proposition implies Theorem~\ref{resomf}. We assume by contradiction that $f_{\ast}(\pi_{1}(X))$ has infinite index in the group $\pi_1(\Sigma_1)\times \pi_1(\Sigma_2)\times \pi_1(\Sigma_3)$ and apply Theorem~\ref{clfsp}. Let $\Sigma'_{i}$, $q_i$, $E$ and $X_0$ be as in that theorem and let $p_{i,0}: X_0 \to \Sigma'_{i}$ be the lift of $X\to \Sigma_i$ with respect to the covering $X_0 \to X$. Let $p=(p_{1,0},p_{2,0},p_{3,0})$ and $Y_0=p(X_0)$. According to Theorem~\ref{clfsp}, there exists a point $o\in E$ such that
$$Y_0=(q_1+q_2+q_3)^{-1}(o).$$
Since the maps $p_{i,0} : X_0\to \Sigma'_i$ are distinct Kodaira fibrations, one can apply Proposition~\ref{normalisation} to $X_0$ instead of $X$. We observe that if we fix a point $m$ of $\Sigma_{1}'\times \Sigma_{2}'\times \Sigma_{3}'$ and choose suitable local coordinates on the $\Sigma_{i}'$, the map $q_{1}+q_{2}+q_{3}$ is given locally near $m$ by:
\begin{equation}\label{formappldddd}
(q_{1}+q_{2}+q_{3})(z_{1},z_{2},z_{3})=z_{1}^{k_{1}}+z_{2}^{k_{2}}+z_{3}^{k_{3}},
\end{equation}
where the $k_i$'s are nonzero natural integers. Its critical points are thus isolated. If the point $o$ is not a regular value of the map $q_1+q_2+q_3$, then Equation~\eqref{formappldddd} shows that $Y_0$ has isolated singular points, a contradiction with Proposition~\ref{normalisation}. Hence $o$ must be a regular value. We then apply Lemma~\ref{lemSmooth} (to $X_0$ instead of $X$) and obtain that the restriction of $p_{1,0}$ to $Y_0$ is a submersion. We now pick a point $(s_1,s_2,s_3)\in X_{0}\subset \Sigma'_1\times \Sigma'_2\times \Sigma'_3$ such that $s_2$ and $s_3$ are critical points of $q_2$ and $q_3$ respectively. This is possible since $(p_{2,0},p_{3,0}) : X_0 \to \Sigma'_2 \times \Sigma'_3$ is surjective. Since $o$ is a regular value, $s_1$ is a regular point of $q_1$. We consider the following commutative diagram:

$$\xymatrix{X_{0} \ar[r] \ar[d] & \Sigma'_{2}\times \Sigma'_{3} \ar[d]^{q_{2}+q_{3}} \\
\Sigma'_{1} \ar[r]^{o-q_{1}} & E. \\ }$$

\noindent The composition of the left and bottom arrow is a submersion at $(s_1,s_2,s_3)$. This gives a contradiction since $(s_2,s_3)$ is a critical point of the right arrow. This concludes the proof of Theorem~\ref{resomf}.

\begin{rem} Here is another argument to rule out the case where $o$ is a regular value of $q_1+q_2+q_3$. Although it is less direct, we find it interesting to note it. If $o$ is a regular value, we have seen that $Y_0$ is a Kodaira fibration, hence it is aspherical. But according to~\cite{bhms2002,dps,Llo-fourier} the fundamental group of a smooth fiber of $q_1+q_2+q_3$ is never ${\rm FP}_{3}$ (see~\cite{browncoho} for the definition of this property), hence cannot be the fundamental group of a closed aspherical manifold. 
\end{rem}

We now turn to the proof of Proposition~\ref{normalisation}. First observe that if $Y$ has an isolated singular point $y_0$, then it must be normal near $y_0$ since $Y$ is a hypersurface locally embedded in $\C^3$ (see e.g. Theorem 3.1 in~\cite{laufer}). Hence $Y$ coincides with its normalization near $y_0$ and the smoothness of $\widetilde{Y}$ implies the second affirmation of the proposition.

We now prove that $\widetilde{Y}$ is smooth. Since the singular set of a normal space has codimension at least $2$ (see~\cite[p. 128]{gr}), $\widetilde{Y}$ has at most isolated singularities. We will thus work near a fixed (but arbitrary) point $\widetilde{y}$ of $\widetilde{Y}$ and apply the material presented in Section~\ref{isosing} to the germ $(\widetilde{Y},\widetilde{y})$ to prove that $\widetilde{Y}$ is smooth at $\widetilde{y}$. 

We fix a lift $\widetilde{f} : X \to \widetilde{Y}$ of the map $f : X \to Y$ with respect to the natural map $\widetilde{Y}\to Y$ and denote by $h_1$ the map that makes the following diagram commutative.
$$
\xymatrix{\widetilde{Y} \ar[r] \ar[drr]^{h_{1}} & Y \ar[r] & \Sigma_{1}\times \Sigma_{2} \times \Sigma_{3} \ar[d] \\
 & & \Sigma_{1} \\}
$$

\noindent Let $\widetilde{y}$ be a point of $\widetilde{Y}$. After fixing a chart centered at $h_1(\widetilde{y})$ in $\Sigma_1$ we think of $h_1$ as a map
$$(\widetilde{Y}, \widetilde{y})\to (\mathbb{C},0).$$   
We shall study its Milnor fiber using the material described in Section~\ref{isosing}. 

\begin{lemma} The boundary of the Milnor fiber of $h_1$ at $\widetilde{y}$ is connected. 
\end{lemma}
{\it Proof.} We fix a germ of function $r : (\widetilde{Y},\widetilde{y})\to [0,+\infty)$ which defines $\widetilde{y}$, as defined at the beginning of Section~\ref{isosing}. We can assume that for a suitable embedding of $(\widetilde{Y},\widetilde{y})$ in $(\C^N, 0)$, $r$ is the square of the Euclidean norm. We observe that according to the first point of Theorem~\ref{thmIsolsing}, the boundary of the Milnor fiber of $h_{1}$ at $\widetilde{y}$ (which is a disjoint union of circles) is diffeomorphic to the set $h_{1}^{-1}(0) \cap \{ r=\varepsilon\}$ for small enough $\varepsilon$. Now $(h_{1}^{-1}(0),\widetilde{y})$ is a germ of complex analytic curve. Hence $h_{1}^{-1}(0) \cap \{ r=\varepsilon\}$ is the intersection of the curve $h_{1}^{-1}(0)$ with a small sphere. Its number of connected components is the number of irreducible components of $h_{1}^{-1}(0)$ at $\widetilde{y}$. 

Hence, the statement of the lemma is equivalent to the fact that the curve $h_{1}^{-1}(0)$ is irreducible at $\widetilde{y}$, a fact that we will now prove. At the end the conclusion will follow from the fact that a holomorphic map $\varphi$ from a disc centered at $0$ in $\C$ to a complex curve has its image contained in one (local) irreducible  component of the curve at $\varphi (0)$. Hence it cannot be open if the target curve is reducible.

The map $\widetilde{f}$ is open, since it is finite-to-one and $\widetilde{Y}$ is irreducible at all of its points (see \S 4 of Chapter 5 in~\cite{gr}). Since $h_1\circ \widetilde{f} = f_1$ this implies that the induced map $f_{1}^{-1}(0)\to \widetilde{Y} \cap h_{1}^{-1}(0)$ is open. We pick a point $x\in \widetilde{f}^{-1}(\widetilde{y})$ and a connected neighbourhood $U$ of $x$ in $f_{1}^{-1}(0)$ such that $U\cap \widetilde{f}^{-1}(\widetilde{y})=\{x\}$. We may assume that $U$ is a disc, because $f^{-1}(0)$ is smooth. Since $h_{1}^{-1}(0)$ has dimension $1$, if it was reducible at $\widetilde{y}$, $\widetilde{f}(U-\{x\})$ and hence $\widetilde{f}(U)$ would thus be contained in one local irreducible component of $h_{1}^{-1}(0)$. This would contradict the fact that the map $f_{1}^{-1}(0)\to \widetilde{Y} \cap h_{1}^{-1}(0)$ induced by $\widetilde{f}$ is open. Hence $h_{1}^{-1}(0)$ is irreducible at $\widetilde{y}$. This implies that the boundary of the Milnor fiber of $h_1$ is a circle.\hfill $\Box$

According to the lemma, the Milnor fiber of $h_1$ at $\widetilde{y}$ is a compact Riemann surface with one boundary component. We shall prove that its genus is zero, hence that it is is a disc. Proposition~\ref{lemMilnDiscSmooth} then implies that $\widetilde{Y}$ is smooth at $\widetilde{y}$, the desired result.

We choose again a point $x\in X$ with $\widetilde{f}(x)=\widetilde{y}$. We fix a chart centered at $x$ in $X$ in such a way that the function $f_1$ locally coincides with the function $(z_1, z_2)\mapsto z_2$ in the chart. Let $\Delta (a)$ be the disc of radius $a>0$ in $\mathbb{C}$. We assume that the polydisc $\Delta(a_0)^2$ is contained in the chart for some $a_0 >0$ and we identify this polydisc with an open set of $X$. We fix a good proper representative $Y_{D, r\le \varepsilon}$ for  $h_1$ (see Section~\ref{isosing}). Since $f$ is open we can assume, up to shrinking $\varepsilon$ and $D$, that $Y_{D, r\le \varepsilon} \subset f(\Delta(a_0)^{2})$. We fix $s\in D-\{0\}$. Define
$$X_{s,r\le \varepsilon}$$
to be a connected component of the preimage under $f|_{\Delta(a_0)^{2}}$ of the Milnor fiber $Y_{s,r\le \varepsilon}$. We can always assume that the map
$$X_{s, r\le \varepsilon}\to Y_{s,r\le \varepsilon}$$
induced by $f$ has no critical point near the boundary of $X_{s,r\le \varepsilon}$. Hence it is a ramified covering. Since $X_{s, r\le \varepsilon}$ is contained in the disc $\Delta(a_0)\times \{s\}$, it has genus zero. This implies that $Y_{s,r\le \varepsilon}$ also has genus zero. Since we have already seen that it has a connected boundary, it must be a disc. This concludes the proof of Proposition~\ref{normalisation}.

\begin{rem}\label{uneseulesubmersion} Let $f_{1} : X \to \Sigma_{1}$ be a Kodaira fibration. Assume that $f_{2} : X \to \Sigma_{2}$ and $f_{3} : X \to \Sigma_{3}$ are fibrations, possibly with critical points and multiple fibers. Assume that $f_1$, $f_2$ and $f_3$ are pairwise distinct. Then the image of $\pi_{1}(X)$ in $\pi_{1}(\Sigma_1)\times \pi_{1}(\Sigma_{2})\times \pi_{1}(\Sigma_{3})$ under the map $(f_1,f_2,f_3)$ has finite index. Indeed, by taking finite coverings, one can first reduce to the case where $f_{2}$ and $f_{3}$ have no multiple fibers. Once this has been done, the proof of Theorem~\ref{resomf} presented above applies verbatim. Indeed during that proof we only used that one of the three fibrations was a submersion. 
\end{rem}


\section{Kodaira fibrations and ${\rm CAT}(0)$ spaces}\label{cat}

This section first discusses a few classical facts concerning curvature and ramified coverings, which allow to build examples of surface-by-surface groups which are ${\rm CAT}(0)$. Afterwards, we prove Theorem~\ref{theo-cat}. 

Consider a Riemannian manifold $M$ of nonpositive curvature and a smooth totally geodesic submanifold $S\subset M$ of codimension $2$. If $X\to M$ is a ramified covering whose branching locus is $S$, one can lift the metric from $M$ to $X$. One obtains a singular tensor which, however, defines a true distance on $X$. Gromov~\cite{gromov} proved that this distance on $X$ is locally ${\rm CAT}(0)$. See also~\cite{charneydavis,pansu} for historical references and further developments around the notions of ramified covering and nonpositive curvature. Now if $\Sigma_1$ and $\Sigma_2$ are two Riemann surfaces of genus greater than $1$, endowed with a hyperbolic metric, and if $f : \Sigma_1 \to \Sigma_2$ is a holomorphic covering map (hence a local isometry), the graph of $f$
$$Graph(f)\subset \Sigma_1 \times \Sigma_2$$
is totally geodesic for the product metric on $\Sigma_1 \times \Sigma_2$. In particular any ramified covering of $\Sigma_1 \times \Sigma_2$, ramified along a disjoint union of such graphs, carries a locally ${\rm CAT}(0)$ metric. More generally, any smooth Riemann surface $D\subset \Sigma_1 \times \Sigma_2$ with the property that both projections $D\to \Sigma_i$ are etale is totally geodesic in $\Sigma_1\times \Sigma_2$. These observations imply that the Atiyah-Kodaira examples as well as all {\it double etale Kodaira fibrations} in the sense of~\cite{cr} carry locally ${\rm CAT}(0)$ metrics. 
This motivates the question: when does a Kodaira fibration carry a locally ${\rm CAT}(0)$ metric? More generally: when is the fundamental group of a Kodaira fibration a ${\rm CAT}(0)$ group? Recall that a group is called ${\rm CAT}(0)$ if it admits a properly discontinuous and cocompact action on a proper ${\rm CAT}(0)$ space.  

\vspace{.4cm}

Theorem~\ref{theo-cat} gives a partial answer to these questions. In the course of its proof, as well as in Sections~\ref{htmcg} and~\ref{geometricpart}, we will repeatedly use the next lemma, whose proof is obvious once one remembers that surface groups have trivial center. 

\begin{lemma}\label{kmetpd} Let $G$ be a group. Assume that $G$ has a normal subgroup $R$ which is a surface group and let $\pi : G \to G/R$ be the quotient morphism. Then the centralizer $\Lambda$ of $R$ in $G$ is normal in $G$. The restriction of $\pi$ to $\Lambda$ is an isomorphism onto the kernel of the monodromy morphism $\varrho$. 
\end{lemma}

We now turn to the proof of Theorem~\ref{theo-cat}. Besides classical results on ${\rm CAT}(0)$ spaces, one of our main tools will be a result by Monod~\cite{monod} concerning actions of direct products on ${\rm CAT}(0)$ spaces. Let $G$ be a surface-by-surface group with infinite monodromy. So we have a short exact sequence 
\begin{equation}\label{extcat}
\xymatrix{0 \ar[r] & R \ar[r] & G \ar[r]^p & Q \ar[r] & 0,}
\end{equation}
where both $R$ and $Q$ are surface groups. We let $\varphi : Q \to \out(R)$ be the natural monodromy morphism. We asssume that $G$ is ${\rm CAT}(0)$ and, by contradiction, that $\varphi$ is not injective. We consider the centralizer of $R$ in $G$, denoted by $\Lambda$. According to Lemma~\ref{kmetpd}, it is a normal subgroup of $G$, the subgroup generated by $R$ and $\Lambda$ is isomorphic to $R\times \Lambda$ and the restriction of $p$ to $\Lambda$ is an isomorphism onto $Ker(\varphi)$. 

\vspace{.4cm}

We fix a properly discontinuous and cocompact action $G\curvearrowright (E,d)$ of $G$ on a proper ${\rm CAT}(0)$ space $(E,d)$. We first assume that the group $R\times \Lambda$ has no fixed point in the visual boundary of $E$. Thanks to Remark 39 in~\cite{monod} we can take a closed $R\times \Lambda$-invariant convex subset $M\subset E$ which is minimal for these properties and {\it canonical}. A careful analysis of the definition of $M$ in \cite{monod} shows that it is $G$-invariant. The $G$-action is still properly discontinuous and cocompact on $M$. We now work with the space $M$. 

Since the action of $R\times \Lambda$ on $M$ is minimal, Corollary 10 in~\cite{monod} implies that there exists an isometric splitting
$$M\simeq M_{1}\times M_{2},$$ 
such that the action of $R\times \Lambda$ is a product action: $R$ acts isometrically on $M_{1}$, $\Lambda$ acts isometrically on $M_{2}$ and the action of $R\times \Lambda$ on $M$ is the product of these two actions. We will now need a slightly more precise result, which follows from Monod's proof. Namely:

\begin{prop} The $G$-action on $M\simeq M_{1}\times M_{2}$ is a product action. 
\end{prop}
{\it Proof.} We explain why this follows from Monod's construction in~\cite[\S 4.3]{monod}. Recall that by definition $M$ is a closed minimal invariant convex set for the $R\times \Lambda$-action. Monod proves that the set of minimal closed $R$-invariant convex subsets of $M$ is non-empty and that the union of all such subsets is closed, convex, $R\times \Lambda$-invariant and splits as a direct product $M_{1}\times M_{2}$ in such a way that the minimal closed convex $R$-invariant subsets are of the form $M_{1}\times \{\ast\}$. He proves moreover that the action of $R\times \Lambda$ is a product action. Since the product action is minimal, the $\Lambda$-action on $M_{2}$ must be minimal. Hence the set of minimal $R$-invariant (resp. $\Lambda$-invariant) closed convex subsets can be identified with $M_{2}$ (resp. $M_{1}$). The identification between $M$ and $M_{1}\times M_{2}$ can be thought of as a map
$$M\to M_{1}\times M_{2}$$ which takes a point $x$ to the pair of minimal closed convex subsets containing $x$ for the respective actions of $R$ and $\Lambda$. Since $G$ normalizes both $R$ and $\Lambda$, it acts naturally on $M_{1}$ and $M_{2}$ and the previous identification shows that the $G$-action on $M$ is a product action.\hfill $\Box$
 
The previous proposition implies that we can now consider the $G$-action on each $M_{i}$ separately. It factors through a faithful action of $G/\Lambda$ on $M_1$ (resp. $G/R$ on $M_2$). 

\begin{prop}\label{pddd} The $G/\Lambda$-action on $M_1$ is properly discontinuous and cocompact. Similarly the $G/R$-action on $M_2$ is properly discontinuous and cocompact.  
\end{prop}
{\it Proof.} The fact that both actions are cocompact is clear. Indeed the $G$-action on $M$ is cocompact. 
Hence the $G$-action on each $M_i$, which are quotients of $M$, is also cocompact.

To show that the actions are properly discontinuous, we make the following observations. The spaces $M_{i}$ are complete locally compact ${\rm CAT}(0)$ spaces. Endowed with the compact open topology their isometry groups are locally compact groups. A subgroup
$$H< {\rm Isom}(X)$$
of the isometry group of a complete locally compact ${\rm CAT}(0)$ space (endowed with the compact open topology) is discrete if and only if its action on $X$ is properly discontinuous (see \cite[5.67]{DruKap-17}). Hence we must show that $G/\Lambda$ and $G/R$ are discrete in the groups ${\rm Isom}(M_1)$ and ${\rm Isom} (M_2)$ respectively. We consider the case of $G/\Lambda$ first. Since $M_{1}$ arises as a convex invariant subset for the action of $R$ on $M$, $R$ acts properly discontinuously on $M_1$. Hence
$$R<{\rm Isom}(M_1)$$
is a discrete subgroup. Thus there exists a neighborhood $U$ of the identity in ${\rm Isom}(M_1)$ such that any element of $U$ which normalizes $R$ must actually centralize $R$. Since the centralizer of $R$ in $G/\Lambda< {\rm Isom}(M_1)$ is trivial, $\left(G/\Lambda\right) \cap U$ is trivial and $G/\Lambda$ is discrete.

The argument is similar for the action of $G/R$ on $M_{2}$. The subgroup $\Lambda < {\rm Isom}(M_2)$ is discrete (as $M_2$ appears as an invariant subset for the action of $\Lambda$ on $M$). Let $\Lambda_0 < \Lambda$ be a 2-generated free subgroup; such a subgroup exists, since $\Lambda$ is a non-trivial normal subgroup of a surface group. Then there is a neighbourhood $U$ of the identity in ${\rm Isom} (M_2)$ such that any element of $U$ which normalizes $\Lambda$ must centralize $\Lambda_0$. Hence, we only have to explain why the centralizer of $\Lambda_0$ in $G/R$ is trivial. But $G/R$ is isomorphic to a surface group, and the centralizer of any non-cyclic subgroup of a surface group is trivial. This concludes the proof.\hfill $\Box$

\vspace{.4cm}

We now fix a point $(m_1,m_2)\in M_{1}\times M_{2}$. Consider the map $f : G \to M_1\times M_2$ defined by $f(g)=(g\cdot m_1,g\cdot m_2)$. This is a quasi-isometry since the action of $G$ on $M$ is properly discontinuous and cocompact. We have a commutative diagram:
$$\xymatrix{G \ar[r]^{f} \ar[d] & M_{1} \times M_{2} \ar[d] \\
G/\Lambda \times G/R \ar[r] & M_{1}\times M_{2} \\}$$
where the vertical map on the right is the identity and the horizontal map on the bottom is given by the orbit maps of $m_{1}$ and $m_2$ for the actions of $G/\Lambda$ and $G/R$. Proposition~\ref{pddd} implies that the bottom map is a quasi-isometry. Hence the vertical arrow on the left must also be a quasi-isometry. But an injective morphism between two finitely generated groups is a quasi-isometry only if its image has finite index. On the other hand the image of the ``diagonal" morphism 
\begin{equation}\label{sedddd}
G \to G/\Lambda \times G/R
\end{equation}
is of infinite index. Indeed, by taking the quotient by the subgroup $R$ on the left and by its image on the right, we obtain a morphism 
$$G/R\to \varphi(Q)\times G/R$$ whose image is the graph of $\varphi$. If the image of the diagonal morphism in~\eqref{sedddd} had finite index, the graph of $\varphi$ would have finite index in $\varphi (Q)\times G/R$, but this only happens if the image of $\varphi$ is finite. We thus obtain a contradiction. 

To conclude the proof, we must now deal with the case when $R\times \Lambda$ fixes a point in the visual boundary of $E$. We will need the following proposition, which is classical (see for instance~\cite{hosaka}). We include the proof for the sake of completeness.

\begin{prop} Let $\Gamma \curvearrowright Z$ be a group acting properly discontinuously and cocompactly on a ${\rm CAT}(0)$ space. Let $N< \Gamma$ be a finitely generated subgroup. If $N$ fixes a point in the visual boundary of $Z$, then the centralizer of $N$ in $\Gamma$ is infinite.  
\end{prop}
{\it Proof.} Let $\{ n_{1}, \ldots , n_{k}\}$ be a finite generating set for $N$. Let $\xi$ be a fixed point for $N$ in the visual boundary of $Z$. We pick a geodesic ray $\alpha : [0,+\infty)\to Z$ representing $\xi$. We pick a constant $C$ such that 
$$d(n_{j} \alpha (t),\alpha (t))\le C$$
for any $t\ge 0$ and any $j\in \{1,\ldots , k\}$. Since the $\Gamma$-action on $Z$ is cocompact, there must exist a constant $A\ge 0$ such that the translates by $\Gamma$ of the ball of radius $A$ centered at $\alpha (0)$ cover $Z$. So we pick a sequence $(\gamma_{i})_{i\ge 0}$ of elements of $\Gamma$ such that $d(\gamma_{i}\alpha (0),\alpha (i))\le A$. We now estimate $d(\gamma_{i}^{-1}n_{s}\gamma_{i} \alpha (0),\alpha (0))$. By the triangle inequality, this is less or equal to:
$$d(\gamma_{i}^{-1}n_{s}\gamma_{i} \alpha (0),\gamma_{i}^{-1}n_{s} \alpha (i))+d(\gamma_{i}^{-1}n_{s} \alpha (i),\gamma_{i}^{-1}\alpha (i))+d(\gamma_{i}^{-1}\alpha (i),\alpha (0)).$$
By our previous choices this is bounded above by $2A+C$. Hence:
$$d(\gamma_{i}^{-1}n_{s}\gamma_{i} \alpha (0),\alpha (0))\le 2A+C.$$
Since the action of $\Gamma$ is properly discontinuous, the set $$\{\gamma_{i}^{-1}n_{s}\gamma_{i}\}_{i\ge 0, 1\le s \le k}$$ must be finite. Hence there exists a subsequence $(\gamma_{i(m)})_{m\ge 0}$ of distinct elements such that for each $s\in \{1, \ldots , k\}$, $\gamma_{i(m)}^{-1}n_{s}\gamma_{i(m)}$ does not depend on $m$. Thus the infinite set $\{\gamma_{i(0)}\gamma_{i(m)}^{-1}\}_{m\ge 1}$ is made of elements commuting with a generating set for $N$; it follows that it is contained in the centralizer of $N$.\hfill $\Box$

\vspace{.4cm}

Let us explain how to conclude from this proposition. Since the monodromy of the extension~\eqref{extcat} is assumed to have infinite image, the kernel of the monodromy morphism is an infinite rank free group. Hence the group $\Lambda$, which is isomorphic to that kernel, is also an infinite rank free group. We fix a basis $B$ for $\Lambda$, pick two distinct elements $b_{1}, b_{2}$ in $B$ and consider the free group $\Lambda_0 < \Lambda$ generated by $b_1$ and $b_2$. If $R\times \Lambda$ fixes a point in the visual boundary of $E$, so does the finitely generated group $R\times \Lambda_0$. We apply the previous proposition to $N=R\times \Lambda_0$ and obtain that the centralizer $C_{G}(N)$ of $N$ is infinite. But this centralizer is the intersection of the centralizer of $R$ in $G$, which is $\Lambda$, with the centralizer $C_{G}(\Lambda_{0})$ of $\Lambda_{0}$ in $G$. We thus get $C_{G}(N)=C_{\Lambda}(\Lambda_{0})$. This is absurd since $C_{\Lambda}(\Lambda_{0})$ is trivial. This concludes the proof of Theorem~\ref{theo-cat}.


\section{Building injective morphisms to the mapping class group}\label{htmcg}

This section contains some preliminary material which will be used in Section~\ref{geometricpart}. In Section~\ref{vaj} we give a criterion to establish the injectivity of certain morphisms to the mapping class group (Proposition~\ref{propositioninjg}) and then introduce {\it polysurface groups} and some of their properties in Section~\ref{ggcepg}. 

\subsection{Extensions by surface groups and their monodromy}\label{vaj}

The ideas presented below are inspired by Johnson's proof of Theorem~\ref{theoremedejohnson}. They will allow us to build new examples of extensions of the form
\begin{equation}\label{extensiondebasebisss}
\xymatrix{0\ar[r] & R \ar[r] & G \ar[r] & G/R \ar[r] & 0}
\end{equation}
 with $R$ a surface group where the monodromy morphism is injective. We thus consider the following algebraic construction, which will reappear as a consequence of a geometric construction in Section \ref{geometricpart}. 
 
 Let $G$ be a group which fits into a short exact sequence

{\small \begin{equation}\label{onetwothreeone}
\xymatrix{0\ar[r] & N \ar[r] & G \ar[r]^{p} & Q \ar[r] & 0,}
\end{equation}}
\noindent where $N$ is a surface-by-surface group with injective monodromy. We assume furthermore that we are given two more exact sequences as follows:

{\small \begin{equation}\label{ttrruuzeroo}
\xymatrix{0 \ar[r] & R_1 \ar[r] & H_{1} \ar[r]^{p_{1}} & Q \ar[r] & 0}
\end{equation}}

{\small \begin{equation}\label{ttrruu}
\xymatrix{0 \ar[r] & R_2 \ar[r] & H_{2} \ar[r]^{p_{2}} & Q \ar[r] & 0}
\end{equation}}

\noindent where the $R_{i}$ are surface groups, and the two extensions have injective mono\-dromy. Furthermore, we assume that we are given surjective morphisms $f_{i} : G\to H_{i}$ and $u_{i} : N\to R_{i}$ for $i=1,2$ such that the kernel of each $u_{i}$ is a surface group, $Ker(u_{1})\neq Ker(u_{2})$ and such that the following diagram is commutative for $i=1,2$:

{\small \begin{equation}
\xymatrix{
0\ar[r] & N \ar[r] \ar[d]^{u_{i}} & G \ar[r]^{p} \ar[d]^{f_{i}} & Q \ar[r] \ar[d]^{id} & 0 \\
0 \ar[r] & R_i \ar[r] & H_{i} \ar[r]^{p_{i}} & Q \ar[r] & 0. \\
}
\end{equation}}

\noindent This implies that the groups $Ker(u_{i})$ are normal in $G$ and that the quotients $G/Ker(u_{i})$ are naturally isomorphic to $H_{i}$. We then have: 
 
\begin{prop}\label{propositioninjg} Under the previous hypotheses, the monodromy of the extension 
{\small \begin{equation}
\xymatrix{0 \ar[r] & Ker(u_{i}) \ar[r] & G \ar[r] & H_{i} \ar[r] & 0}
\end{equation}}
is injective for $i=1,2$. 
\end{prop}

In the course of the proof below we will use that a nontrivial, finitely generated normal subgroup of a surface group must have finite index. Since $Ker(u_{1})$ and $Ker(u_{2})$ are distinct, this implies that $Ker (u_{1})$ cannot be contained in $Ker(u_{2})$ (and vice versa). Indeed if $Ker(u_{1}) <Ker (u_2)$, then the group $Ker (u_2)/Ker (u_1)$ must be finite. Since $R_1$ is torsion-free, we must then have $Ker(u_2)=Ker(u_1)$, and this contradicts our initial hypothesis. We will also need the following:

\begin{lemma}\label{lemmef} Let $R$ be a surface group and let $L$ be a nontrivial normal subgroup of $R$. Let $f$ be an automorphism of $R$. If $f$ is the identity on $L$ then $f$ is the identity on all of $R$. 
\end{lemma}
{\it Proof.} The group $R$ and the automorphism $f$ act on the Gromov boundary $\partial (R)$ of $R$, which is a circle. The map $f$ fixes the limit set of $L$ in $\partial(R)$ pointwise. Since $L\lhd R$ is normal, this limit set is invariant under the action of $R$. By the minimality of the action of $R$ on $\partial (R)$ \cite[Ch. 8, Obs. 27]{GhydlH-90}, it must coincide with $\partial (R)$. Hence $f$ is the identity on $\partial (R)$. This implies that $f$ is the identity on $R$.\hfill $\Box$

\vspace{.4cm}

\noindent {\it Proof of Proposition~\ref{propositioninjg}.} We consider the case $i=1$. We assume by contradiction that the monodromy is not injective. Then $G$ contains the direct product $Ker(u_{1})\times \Lambda_{1}$, where $\Lambda_{1}$ is the centralizer of $Ker(u_{1})$ in $G$, which is nontrivial by Lemma~\ref{kmetpd}. The groups $u_{2}(Ker (u_{1}))$ and $f_{2}(\Lambda_{1})$ commute in $H_{2}$. Since $Ker(u_{1})$ is normal in $N$, $u_{2}(Ker (u_{1}))$ is normal in $R_{2}$. It is also non-trivial since $Ker(u_1)$ is not contained in $Ker(u_2)$. Lemma~\ref{lemmef} then implies that $f_{2}(\Lambda_{1})$ centralizes all of $R_{2}$. Since the monodromy of the extension~\eqref{ttrruu} is injective we must have $f_{2}(\Lambda_{1})=\{e\}$. This implies that $\Lambda_{1}$ and hence $Ker( u_1)\times \Lambda_1$ is contained in $N$. This contradicts the injectivity of the monodromy of $u_1: N\to R_{1}$, thus completing the proof.\hfill $\Box$

\vspace{.4cm}

\subsection{Central extensions of polysurface groups}\label{ggcepg}

Polysurface groups are the groups obtained via iterated extensions by surface groups. More precisely, we define them as follows.

\begin{defi}\label{polys} A group $G$ is a polysurface group of length $n$ if there exists a filtration $(G_{i})_{0\le i \le n}$ with $\{e\}=G_{0}<G_{1}<\ldots < G_{n}=G$ such that $G_i$ is a normal subgroup of $G$ and $G_{i}/G_{i-1}$ is isomorphic to a surface group for each $1\le i \le n$.   
\end{defi}

When $n=2$, we recover the definition of surface-by-surface groups. Note that if $G$ and $(G_{i})_{0\le i \le n}$ are as in the definition, then $G/G_{1}$ and $G_{n-1}$ are polysurface groups of length $n-1$. 

When $n\ge 3$, this definition differs slightly from the one given in~\cite{j99} where the $G_i$'s are only assumed to be normal in $G_{i+1}$\footnote{However, one can prove that a group satisfying Johnson's definition must virtually satisfy our definition. This follows from Johnson's result saying that a given group has finitely many polysurface group structures~\cite{j99}.}. But we will work with the above definition.

We can now state the following proposition, that we will use repeatedly in Section~\ref{atkf}.

\begin{prop}\label{ceatofis} Let $G$ be a polysurface group and let $A$ be a finite abelian group. Then any central extension
$$\xymatrix{0 \ar[r] & A \ar[r] & \Gamma \ar[r] & G \ar[r] & 0}$$
becomes trivial after passing to a finite index subgroup of $G$.
\end{prop}
{\it Proof.} We first prove the assertion if $G$ is a surface group. The extension we consider corresponds to a class in the group $H^{2}(G,A)$. Since $H_{1}(G,\Z)$ is torsion-free, the universal coefficient theorem implies that
$$H^2(G,A)\simeq Hom(H_{2}(G,\Z),A).$$
If $G_{1}<G$ is a subgroup of finite index such that the image of the natural map $H_{2}(G_{1},\Z)\to H_{2}(G,\Z)$ is generated by an element which is divisible in $H_{2}(G,\Z)$ by the order of $A$, then the extension is trivial when restricted to $G_{1}$. This proves the proposition in the case of a surface group.

We now assume that $G$ is a polysurface group of length $n$ and prove the assertion by induction on $n$. We can assume that $n\ge 2$. Let $G_{1}\lhd G$ be a normal subgroup such that $G/G_{1}$ is a polysurface group of length $n-1$. We pick $G_{1}'$ of finite index in $G_1$ such that the extension is trivial over $G_{1}'$. Let $i : G_{1}' \to \Gamma$ be a lift of $G_{1}'$. Let $\Gamma'$ be the normalizer of $i(G_{1}')$ in $\Gamma$. It has finite index in $\Gamma$. So we have an extension 
\begin{equation}\label{hdeuxit}
\xymatrix{0\ar[r] & A \ar[r] & \Gamma' \ar[r] & G' \ar[r] & 0,}
\end{equation}
where $G'$ is the image of $\Gamma'$ in $G$ and the subgroups $A$ and $i(G_{1}')$ are normal in $\Gamma'$. By taking the quotient by $i(G_{1}')$ in~\eqref{hdeuxit}, we obtain a central extension of $G'/G_{1}'$ by $A$:
\begin{equation}\label{hdeuxitt}
\xymatrix{0 \ar[r] & A \ar[r] & \Gamma'/i(G_{1}') \ar[r] & G'/G_{1}' \ar[r] & 0.}
\end{equation}
The group $G'/G_{1}'$ is an extension of a polysurface group of length $n-1$ by a finite group. Indeed it contains the finite group $\left(G_1\cap G'\right)/G_{1}'$, and the quotient is a finite index subgroup of a polysurface group of length $n-1$, hence is a polysurface group of the same length. In particular, $G'/G_{1}'$ admits a finite index subgroup which is a central extension of a finite abelian group by a polysurface group of length $n-1$. Using the induction hypothesis twice, we obtain a finite index subgroup $M<\Gamma'/i(G_{1}')$, containing $A$, on which the extension is trivial. This implies that there is a left-splitting $M\to A$ for the inclusion $A\hookrightarrow M$. It induces a left-splitting $\Gamma''\to A$ of the inverse image $\Gamma''$ of $M$ in $\Gamma'$. Thus, $\Gamma''< \Gamma$ is a finite index subgroup on which the extension is trivial.\hfill $\Box$

\begin{rem} The above proposition is not true for arbitrary groups $G$. Indeed, if $G$ is residually finite and if $\Gamma$ is a {\it central} extension of $G$ by a finite group, then the extension is trivial on a finite index subgroup of $G$ if and only if $\Gamma$ is residually finite. However, there are examples showing that this condition need not be satisfied (see~\cite{deligne}, as well as~\cite[Ch. 8]{ABCKT-95} for more examples). 
\end{rem}


\section{Families of Riemann surfaces with injective monodromy}\label{geometricpart}

In this section we provide two geometric constructions of surface bundles to which we apply the group theoretic results from Section \ref{htmcg}. Roughly, they both rely on performing the classical construction of Kodaira fibrations from~\cite{Ati-69,Kod-67} ``in family". The first construction leads to families of Riemann surfaces with injective monodromy but only applies in the context of iterated Kodaira fibrations. The second construction applies more generally although it produces families with fibers of higher genus.

If $\pi : Z\to B$ is a fiber bundle whose fiber $F$ is a closed oriented surface of genus greater than $1$, then $\pi$ always induces a short exact sequence between the fundamental groups of $F$, $Z$ and $B$. This is due to the fact that the map $\pi_{2}(B)\to \pi_{1}(F)$ coming from the long exact sequence in homotopy induced by $\pi$ is always trivial. Indeed, $\pi_{1}(F)$ does not contain any nontrivial Abelian normal subgroup. The same conclusion applies as long as $\pi_{1}(F)$ has this last property (e.g. if $F$ is a Kodaira fibration). We will regularly use this fact without further reference.

\subsection{Iterated fibrations with a fiberwise involution}\label{atkf}

Sections~\ref{atkf} and~\ref{candfi} are devoted to a geometric construction which will lead to the proof of Theorem~\ref{injm}. The construction is made by induction. So for now, we assume that we are given an iterated Kodaira fibration $X$ of dimension $n$ with injective monodromy. Let 
\begin{equation}\label{offti}
\pi : X \to Y
\end{equation}
be the corresponding fibration, where $Y$ is an iterated Kodaira fibration of dimension $n-1$. In this section we explain how, by taking finite coverings, we can produce another fibration with the same property but which moreover carries a fixed point free involution which preserves the map $\pi$. 

We write $R$ for the kernel of the map
$$(p_{1})_{\ast} : \pi_{1}(X)\to \pi_{1}(Y).$$
Let $R_{2}$ be the kernel of the natural map $R \to H_{1}(R, \Z/2\Z)$. Since $R_2 < R$ is characteristic, we have an extension 
\begin{equation}\label{extuusstt}
\xymatrix{0 \ar[r] & H_{1}(R , \Z/2\Z) \ar[r] & \pi_{1}(X)/R_{2} \ar[r] & \pi_{1}(Y)\ar[r] & 0.}
\end{equation}
After passing to a finite index subgroup of $\pi_{1}(Y)$, we can assume that this extension is central and applying Proposition~\ref{ceatofis} and passing to another finite index subgroup, we can assume that it is actually trivial. So let $Y'\to Y$ be a finite covering such that the above extension is trivial on $\pi_{1}(Y')$. Let $X'\to X$ be the induced covering. We have a surjective morphism $\pi_{1}(X')\to H_{1}(R , \Z/2\Z)$ whose restriction to $R$ is the natural map $R\to H_{1}(R,\mathbb{Z}/2\mathbb{Z})$. Thus we can choose a surjective morphism $\varphi : \pi_1(X')\to \Z/2\Z$ which is nontrivial on $R$. Let $X''$ be the covering space corresponding to the kernel of the morphism $\varphi$. Let $\pi'' : X''\to Y'$ be the composition of the covering map $X''\to X'$ with a lift $\pi': X'\to Y'$ of $\pi$.

\begin{prop} The monodromy of the fibration $\pi '' : X'' \to Y'$ is injective. 
\end{prop}
{\it Proof.} Note that since $\varphi (\pi_{1}(X'))=\varphi (R)$, the map $X''\to Y'$ is $\pi_{1}$-surjective (equivalently, has connected fibers). We now use repeatedly Lemma~\ref{kmetpd}. If the monodromy of the fibration $X''\to Y'$ is not injective, the centralizer $\Lambda$ of $R \cap \pi_{1}(X'')$ in $\pi_{1}(X'')$ is nontrivial. By Lemma~\ref{lemmef}, $\Lambda$ must centralize all of $R$. This contradicts the fact that the monodromy of the original fibration~\eqref{offti} is injective.\hfill $\Box$

\vspace{.4cm}

The group $\Z/2\Z$ acts as a Galois group on $X''$ and this action preserves $\pi''$. Its generator $1\in \Z/2\Z$ thus defines a fixed point free involution of $X''$ which leaves $\pi''$ invariant. We finally observe that the covering map $X''\to X'$ restricts to a covering of degree two  $\pi''^{-1}(y)\to \pi'^{-1}(y)$ for $y\in Y'$, while the genus of the fibers of $\pi'$ is identical to the genus of the fibers of $\pi$. Summing up all pieces, we have proved the following proposition. 

\begin{prop}\label{reductionsteppppp} 
Let $\pi : X\to Y$ be an $n$-dimensional iterated Kodaira fibration with injective monodromy and with fibers of genus $g$. Then there exists a finite covering space $X''\to X$ with the following properties:
\begin{enumerate}
 \item there is an induced submersion $\pi'':X''\to Y'$ equipping $X''$ with the structure of an $n$-dimensional iterated Kodaira fibration with injective monodromy; 
 \item $X''$ carries a fixed point free involution $\sigma: X'' \to X''$ such that $\pi''\circ \sigma =\pi''$;
 \item the fibers of $\pi''$ have genus $1+2(g-1)$.
\end{enumerate}
\end{prop}

\subsection{Constructing an $(n+1)$-dimensional fibration}\label{candfi} We now proceed with the induction step. We assume that there exists an $n$-dimensional iterated Kodaira fibration with injective monodromy and with fibers of genus $g$ (for some $n\ge 2$). Thanks to Proposition~\ref{reductionsteppppp}, we can pick an iterated Kodaira fibration $X$ of dimension $n$
$$\pi : X \to Y$$
with injective monodromy and with fibers of genus $1+2(g-1)$, and assume that $X$ is endowed with a holomorphic fixed point free involution $\sigma : X \to X$ such that $\pi \circ \sigma =\pi$. We now repeat the arguments from the previous section. We call again $R < \pi_{1}(X)$ the fundamental group of the fiber and $R_{2}$ the kernel of the map $R \to H_{1}(R , \Z/2\Z)$. First one can pick a finite covering $Y'\to Y$ and the induced covering $X'\to X$ such that the extension~\eqref{extuusstt} is trivial. We fix a lift $X' \to Y'$ of $\pi$, called $\pi'$; observe that the fibers of $\pi'$ also have genus $1+2(g-1)$. We fix a morphism $\phi : \pi_{1}(X')\to H_{1}(R , \Z/2\Z)$ which induces the canonical morphism $R\to H_{1}(R , \Z/2\Z)$ on $R <\pi_1(X')$. Let $X''$ be the covering space corresponding to the kernel of $\phi$. Let $\sigma ' : X' \to X'$ be the lift of $\sigma$ which preserves the map $\pi'$. We call $f : X'' \to X'$ the covering map. 
Let $D$ be the union of the graph of $f$ and of the graph of $\sigma ' \circ f$. Note that $D$ naturally sits as a smooth divisor in the fiber product
$$Z:= X''\times_{Y'}X'=\{(x,z)\in X''\times X', \pi '(f(x))=\pi '(z)\}.$$
We denote by $Z_{y}$ the fiber of $Z$ above a point $y\in Y'$. It is the direct product of the fibers $X''_y$ of $X''\to Y'$ and $X'_y$ of $X'\to Y'$. Note that by construction $f$ restricts to the covering $X''_y\to X'_y$ corresponding to the kernel of the morphism $\pi_1(X'_y)\to H_1(X'_y,\Z/2\Z)$. In particular, this means that we can perform the Atiyah-Kodaira construction on $Z_y$. In the following, we will simply perform this construction ``in family". 

\begin{prop}\label{homologyzeromodtwo} After possibly passing to a finite covering space $Y'' \to Y'$ and taking the induced covering spaces $Z'\to Z$ and $D'\to D$ we can assume that there exists a morphism 
$$\pi_{1}(Z'-D')\to \Z/2\Z$$
which is nontrivial on any small transversal loop enclosing a connected component of $D'$. 
\end{prop}

Observe that given a point $y\in Y$, we can always build a morphism 
\begin{equation}\label{eqnKodFib}
 \pi_{1}(Z_{y}-Z_{y}\cap D)\to \Z/2\Z
\end{equation}
 which is nontrivial on a small loop enclosing a component of $Z_{y}\cap D$. This follows from the fact that the homology class of $Z_{y}\cap D$ is divisible by $2$ in $H_{2}(Z_{y},\Z)$, see~\cite{Ati-69, Kod-67} as well as~\cite{chen} for a detailed discussion. This induces a ramified covering $W_y\to Z_y$ of degree $2$, with ramification locus $Z_y\cap D$, together with two submersions $W_y\to X''_y$ and $W_y\to X'_y$, equipping $W_y$ with the structure of a double Kodaira fibration. Note that the fibers of $W_y\to X''_y$ are ramified coverings of degree $2$ of $X'_y$ whose ramification locus consists of two points; in particular their genus is $2+4(g-1)$.

The point of Proposition~\ref{homologyzeromodtwo} is that the morphism constructed by Kodaira can be virtually extended to $\pi_{1}(Z-D)$. 

\vspace{.4cm}

\noindent {\it Proof of Proposition~\ref{homologyzeromodtwo}.} This is similar to the arguments presented in Section~\ref{atkf} except that we consider the bundle $Z-D\to Y'$. We denote by $L$ the fundamental group of the fiber of this new bundle, by $L_2$ the kernel of the natural map $L\to H_{1}(L,\Z/2\Z)$ and substitute the extension~\eqref{extuusstt} by the following one:
\begin{equation}\label{extuusstttt}
\xymatrix{0 \ar[r] & H_{1}(L,\Z/2\Z) \ar[r] & \pi_{1}(Z-D)/L_{2} \ar[r] & \pi_{1}(Y') \ar[r] & 0.}
\end{equation}
By Proposition \ref{ceatofis} one can take a finite covering space $Y''\to Y'$ such that the restriction of the extension~\eqref{extuusstttt} to $\pi_{1}(Y')$ is trivial. Let then $Z'\to Z$ be the induced covering space and let $D'$ be the inverse image of $D$ in $Z'$. Note that $D’$ is a disjoint union of two connected components which cover the two connected components of $D$. By construction there exists a surjection
$$\pi_{1}(Z'-D')\to H_{1}(L,\Z/2\Z)$$
which splits the extension~\eqref{extuusstttt} restricted to $\pi_{1}(Z'-D')$. Composing it with a morphism from $H_{1}(L,\Z/2\Z)$ to $\Z/2\Z$ which behaves as in \eqref{eqnKodFib} on the fibers of $Z' - D' \to Y''$ gives the desired morphism.\hfill $\Box$

\vspace{.4cm}

We now fix a morphism $\varphi : \pi_{1}(Z'-D')\to \Z/2\Z$ as in the proposition. Its kernel defines a double covering of $Z'-D'$ which extends to give a double ramified covering $Z^{\ast}\to Z'$. The map $Z^{\ast} \to Z'\to Y''$ is $\pi_{1}$-surjective; we denote it by $p$. Observe that by construction the fibers of $Z'-D'\to Y''$ are diffeomorphic to the fibers of $Z-D\to Y'$. We thus have a holomorphic submersion
$$p : Z^{\ast} \to Y''$$
whose fibers are double Kodaira fibrations $Z^*_y$, which are diffeomorphic to the double Kodaira fibration $W_y$ defined above. In particular the fibers of $p$ are Kodaira fibrations with injective monodromy and one of the submersions has fibers of genus $2+4(g-1)$. 

We are now in the position to apply the results from Section~\ref{vaj}. The short exact sequence
\begin{equation}\label{onetwothreeoneconcrete}
\xymatrix{0\ar[r] & \pi_{1}(Z^{\ast}_{y}) \ar[r] & \pi_{1}(Z^{\ast}) \ar[r] & \pi_{1}(Y'')\ar[r] & 0}
\end{equation}
will play the role of the Sequence~\eqref{onetwothreeone} from Section~\ref{vaj}. Let $q_1: X_{1}\to Y''$ (resp. $q_2: X_{2}\to Y''$) be the fibration obtained from the fibration $X'\to Y'$ (resp. $X''\to Y'$) by the base change $Y''\to Y'$. The short exact sequences on fundamental groups built from these two new fibrations play the role of Sequences~\eqref{ttrruuzeroo} and~\eqref{ttrruu}. The corresponding monodromy morphisms are injective since the fibration $\pi : X \to Y$ we started with at the beginning of this section has injective monodromy. We finally have to specify the morphisms $f_1$, $f_2$, $u_1$ and $u_2$. Note that $Z'$ is naturally identified with the fiber product
$$X_{1} \times_{Y''}X_{2}.$$

\noindent The morphism $f_i$ is simply the map on fundamental groups induced by $Z^{\ast}\to Z'\to X_{i}$, $u_{i}$ being its restriction to the fundamental group of the fiber of $p$. The fiber of the projection $Z'\to X_1$ over $x_1\in X_1$ is isomorphic to the fiber $X_{2,q_1(x_1)}$ of $q_2$ over $q_1(x_1)\in Y''$, which is connected. In particular, the fiber of $Z^{\ast}\to X_1$ is identical with the fiber of the Kodaira fibration $Z^{\ast}_{q_1(x_1)}\to X_{1,q_1(x_1)}$ and thus also connected. It follows that $f_1$ is $\pi_1$-surjective and an analogous argument shows that the same holds for $f_2$. Moreover, the projections $u_i$ are precisely the surjective morphisms corresponding to the two Kodaira fiberings of $Z^{\ast}_{q_1(x_1)}$. Proposition~\ref{propositioninjg} now implies that the monodromy of the fibrations $Z^{\ast}\to X_{i}$ is injective for $i=1, 2$. This concludes the proof of Theorem~\ref{injm}. 

\vspace{.4cm}

Finally, let us make the following observations about our proof. We started in Section~\ref{atkf} with any $n$-dimensional iterated Kodaira fibration with injective monodromy and fibers of genus $g$. We proved there that, up to changing $X$ and $Y$ by finite covers, one could assume that $X$ carried a fiberwise holomorphic involution. In particular the fibration $X$ we started with at the beginning of Section~\ref{candfi} can be taken to be a finite covering space of any prescribed iterated Kodaira fibration of dimension $n$ with injective monodromy. The bases of the fibrations $Z^{\ast}\to X_i$ we end up with being again finite covering spaces of $X$ and the fibration $Z^*\to X_2$ having fibers of genus $2+4(g-1)$, we have proved:  

\begin{main}\label{allve} Let $X$ be an $n$-dimensional iterated Kodaira fibration with injective monodromy and fibers of genus $g$. Then there exists a finite covering space $X_0$ of $X$ and a family $Z\to X_0$ of closed Riemann surfaces of genus $2+4(g-1)$ whose monodromy is injective. In particular the fundamental group of $X$ virtually embeds into the mapping class group of a closed surface of genus $2+4(g-1)$. 
\end{main}

Note that this statement of course contains Theorem~\ref{injm2}. It would be of interest to make Theorem~\ref{allve} more effective. Given $X$ as above, what are the minimal degree of the covering $X_{0}\to X$ and the minimal genus of the fibers of the family $Z\to X_{0}$? It would be interesting to answer this question even for specific examples of iterated Kodaira fibrations $X$.   

\begin{rem}
The construction described in this section can also be performed in the category of real smooth manifolds. This might allow to construct more real manifolds which are smooth surface bundles with injective monodromy and whose base is an iterated surface bundle. 
\end{rem}

\subsection{Virtual representations of the mapping class group}\label{sectionmcg}

The construction presented in the previous sections used that the base of all the families of Riemann surfaces that we considered was an iterated Kodaira fibration at only one place: through the use of Proposition~\ref{ceatofis}. Here we will deal with arbitrary bases, by generalizing an argument due to Miller~\cite[\S 3]{miller1984}. The only price to pay is that the families that we produce have larger genus. More precisely, we will prove the following. 
  
\begin{main}\label{basearbitraire} Let $\pi : Z\to B$ be a holomorphic family of closed Riemann surfaces of genus $g\ge 2$ over a complex manifold $B$. Assume that $\pi$ has injective monodromy. Then there exists a finite covering $Z'\to Z$ and a holomorphic family $p : W \to Z'$ of closed Riemann surfaces of genus $2+8(g-1)$ with injective monodromy. Moreover the monodromy of $p$ restricted to any finite index subgroup of $\pi_{1}(Z')$ is irreducible. 
\end{main}

Let $\pi : Z\to B$ be as in the statement of the theorem. Let as before $R \lhd \pi_{1}(Z)$ be the fundamental group of a fiber of $\pi$. We write $[R,R]$ for the kernel of the natural map $R \to H_{1}(R, \mathbb{Z})$ and, for an integer $m$, $R_{m}$ for the kernel of the map $R \to H_{1}(R, \mathbb{Z}/m\mathbb{Z})$. Although the extension
\begin{equation}\label{extensionespacetotalriemann}
\xymatrix{0\ar[r] & R \ar[r] & \pi_{1}(Z)\ar[r] & \pi_{1}(B)\ar[r] & 0}
\end{equation}
need not have a section, we first observe that the corresponding extension
\begin{equation}\label{extensionh1}
\xymatrix{0\ar[r] & H_{1}(R, \mathbb{Z}) \ar[r] & \pi_{1}(Z)/[R,R] \ar[r] & \pi_{1}(B)\ar[r] & 0}
\end{equation}
virtually has a section. This is the content of the next proposition. This is well-known to experts, but we include a proof based on Morita's work~\cite{morita1}. 

\begin{prop}\label{jacobiii} There exists a finite covering space $B_0\to B$ such that the extension~\eqref{extensionh1} has a section over $\pi_{1}(B_{0})$. 
\end{prop}

Note that Earle~\cite{earle} has proved that one can build a family $E\to B$ of complex tori and an embedding $Z\hookrightarrow E$ preserving the projection onto $B$, which coincides (up to translation) with the Jaocbi map in each fiber. The bundle $E\to B$ need not have a continuous section in general. Indeed, Earle describes precisely when such a section exists, see~\cite[\S 8.1]{earle}. Morita has studied the same kind of questions when one starts with a surface bundle which need not have a complex structure, see~\cite{morita1,morita2}. Here we explain how to deduce Proposition~\ref{jacobiii} from their work. We also refer the reader to~\cite{claudon} for a similar discussion for general families of complex tori (not necessarily coming from a family of Riemann surfaces). 

\noindent {\it Proof.} We will use the fact that there exists a map $u : \pi_{1}(Z) \to H_{1}(R,\mathbb{Z})$ which satisfies the following properties:
\begin{enumerate}
\item $u$ is a cocycle i.e. $u(h_{1}h_{2})=u(h_{1})+(h_{1})_{\ast} u(h_{2})$, where $(h_{1})_{\ast}$ is the automorphism of $H_{1}(R,\mathbb{Z})$ induced by $h_1$; 
\item $u(n)=(2g-2)[n]$ for $n\in R$, where $[n]$ stands for the homology class of $n$.
\end{enumerate}
To prove that such a map $u$ exists it is enough to prove it in the case where one considers the so-called Birman exact sequence
$${\small \xymatrix{0\ar[r] & \pi_{1}(S_{g}) \ar[r] & \mcg (S_{g},\ast) \ar[r] & \mcg (S_{g}) \ar[r] & 0.}}$$
Indeed the exact sequence~\eqref{extensionespacetotalriemann} is induced by the one above through a homomorphism from $\pi_{1}(Z)$ to $\mcg (S_{g},\ast)$. The fact that such a map exists for the Birman exact sequence follows from Morita's work, see~\cite[\S 6]{morita1} or ~\cite[\S 1]{morita2}. Another construction is also contained in Earle's work~\cite{earle}. 

We fix a cocycle $u$ satisfying the above conditions and call $\phi : \pi_{1}(Z)/[R,R]\to \pi_{1}(B)$ the morphism appearing in~\eqref{extensionh1}. The cocycle relation, together with the fact that $u(n)=0$ if $n\in [R,R]$, implies that $u$ descends to a cocycle $\underline{u} : \pi_{1}(Z)/[R,R]\to H_{1}(R,\mathbb{Z})$. Let $$N:=\underline{u}^{-1}((2g-2)H_{1}(R,\mathbb{Z}))$$ and let $N_0:=\underline{u}^{-1}(\{0\})$. Observe that $N$ and $N_0$ are subgroups of $\pi_{1}(Z)/[R,R]$ and $N_0<N$. The group $N$ has finite index in $\pi_{1}(Z)/[R,R]$. For each $h\in N$, there exists $c\in H_{1}(R,\mathbb{Z})<N$ such that $\underline{u}(h)=\underline{u}(c)$. Hence $c^{-1}h\in N_0$. This implies that $N$ and $N_0$ have the same image in $\pi_{1}(B)$. Since $N_0$ does not intersect the subgroup $H_{1}(R,\mathbb{Z})\lhd \pi_{1}(Z)/[R,R]$, $\phi$ induces an isomorphism between $N_{0}$ and $\phi (N)=\phi (N_0)$. Choosing $B_0\to B$ to be the covering space corresponding to the subgroup 
$\phi (N)<\pi_{1}(B)$ concludes the proof.\hfill $\Box$

In view of the previous proposition, we now assume that the base $B$ has been replaced by a finite cover (again denoted by the same letter) so that the extension~\eqref{extensionh1} has a section. In particular any of the extensions
\begin{equation}
\xymatrix{0\ar[r] & H_{1}(R, \mathbb{Z}/m\mathbb{Z}) \ar[r] & \pi_{1}(Z)/R_{m} \ar[r] & \pi_{1}(B)\ar[r] & 0}
\end{equation}
has a section. We now let $m=2$. Let $B'\to B$ be the finite covering space corresponding to the kernel of the representation
$$\pi_{1}(B)\to {\rm Aut}(H_{1}(R, \mathbb{Z}/2\mathbb{Z}))$$
and let $Z'\to Z$ be the induced covering. The extension 
\begin{equation}
\xymatrix{0\ar[r] & H_{1}(R, \mathbb{Z}/2\mathbb{Z}) \ar[r] & \pi_{1}(Z')/R_{2} \ar[r] & \pi_{1}(B')\ar[r] & 0}
\end{equation}
is a central extension with a section hence is trivial. Using this and arguing exactly as in Section~\ref{atkf}, one obtains:

\begin{center} {\it Up to replacing $Z$ and $B$ by finite covering spaces and up to replacing $g$ by $k=1+4(g-1)$ we can and do assume that there exists a free $(\mathbb{Z}/2\mathbb{Z})^{2}$-action on $Z$ which preserves the fibers of $\pi$.} 
\end{center}

Note that the only difference compared to Section~\ref{atkf} is that to establish the above fact, we did not make use of Proposition~\ref{ceatofis} which is not available for an arbitrary base $B$. Instead we used that the extension~\eqref{extensionh1} has a section. Note also that we will need below to have two commuting involutions acting on $Z$ (instead of just one involution as in~\ref{atkf}); this trick is due to Miller~\cite[\S 4]{miller1984}. These are provided by the $(\Z/2\Z)^2$-action on $Z$; we will call them $a_1$ and $a_2$.

We now consider the fiber product 
$$W:=Z \times_{B} Z.$$
 Let $D\subset W$ be the union of the graph of the identity map ${\rm id}_Z :Z\to Z$ and the graph of $a_{1}$. This is a smooth divisor in $W$. As before we want to construct a suitable morphism $\pi_{1}(W-D)\to \mathbb{Z}/2\mathbb{Z}$, possibly after changing $W$ and $D$ by finite covering spaces. This is the analogue of Proposition~\ref{homologyzeromodtwo}. In the following proposition, we think of $W$ as a bundle over $Z$, via the first projection $W\to Z$. Consider the bundle $W - D \to Z$ and let $F$ be its fiber. This is a twice-punctured surface.

\begin{prop}\label{prop:classdivisibleby2} There exists a finite covering $f: Z'\to Z$ such that the finite covering $h : W^1\to W$ induced by base change has the property that for $D^1:=h^{-1}(D)$ there is a morphism $\pi_1(W^{1}-D^{1})\to \Z / 2 \Z$ which is nontrivial on any small simple loop enclosing a component of $D^1$.
\end{prop}

\noindent {\it Proof.} The key point is that the bundle $W - D \to Z$ has a section; this is why we took a covering with group $(\mathbb{Z}/2\mathbb{Z})^{2}$ at the beginning of this section. A section is given by the graph of the map $a_{2}$; the fact that this graph does not intersect $D$ comes from the fact that $a_1$ and $a_2$ generate a free $(\mathbb{Z}/2\mathbb{Z})^{2}$-action. In particular the corresponding extension of fundamental groups is a semidirect product. This implies as before that if $f: Z'\to Z$ is the covering corresponding to the kernel of the action of $\pi_{1}(Z)$ on the first homology group of the fiber $F$ with $\mathbb{Z}/2\mathbb{Z}$ coefficients, and if $W^{1}\to W$ is the induced covering, then the group $\pi_{1}(W^{1}-D^{1})$ surjects onto $H_{1}(F,\Z/2\Z)$ in such a way that the restriction of this morphism to the fundamental group of $F$ is the natural map $\pi_{1}(F)\to H_{1}(F,\Z/2\Z)$. Since a small simple loop enclosing a component of $D$ and contained in a fiber of $F$ is nontrivial in the first homology group of $F$ this gives the desired result.\hfill $\Box$ \vspace{.5cm}

Let $W^{2}$ be a double covering space of $W^{1}$ ramified along $D^1$ (such a covering exists by Proposition \ref{prop:classdivisibleby2}). 

We will now explain how to identify $W^{1}$ with a suitable fiber product and $D^1$ with a union of two graphs in this fiber product. This will put us in a very similar situation as in Section \ref{candfi} and thus enable us to apply Proposition~\ref{propositioninjg} to the composition $W^2\to W^1\to Z'$. 

Let $B'\to B$ be the finite covering corresponding to the subgroup $\pi_{\ast}(\pi_1(Z')) < \pi_1(B)$ and let $Z''\to Z$ be the covering obtained by base change with respect to $B'\to B$. The involution $a_1$ lifts to an involution $a_1''$ of $Z''$ which preserves the fibers of $Z''\to B'$. These fibers are homeomorphic to the fibers of $Z\to B$ and thus have genus $1+4(g-1)$. There is a covering $f':Z'\to Z''$ that commutes with the projection to $B'$ and restricts to a covering $Z'_b\to Z''_b$ between the fibers $Z'_b$ of $Z'\to B'$ and $Z''_b$ of $Z''\to B'$. It follows that there is a natural identification of $W^{1}$ with the fiber product 
$$Z'\times_{B'} Z'',$$ 
and the fiber of the projection $W^{1}\to B'$ above $b\in B$ is thus a product $W^{1}_b=Z'_b\times Z''_b$ of (connected) closed Riemann surfaces. Moreover, this identification maps $D^{1}$ to the union of the graphs of $f'$ and $a_1''\circ f'$.

The family $W^1-D^1\to Z'$ restricts to a family $W^1_b - \left(W^{1}_b\cap D^{1}\right)\to Z'_b$. Thus, the ramified double covering $W^{2}\to W^{1}$ restricts to a double covering $W^{2}_b\to W^{1}_b=Z'_b\times Z''_b$ ramified in the smooth divisor $W^{1}_b \cap D^{1}$. This equips the fibers $W^{2}_b$ of the projection $W^{2}\to B'$ with the structure of double Kodaira fibrations; in particular their two projections onto $Z'_b$ and $Z''_b$ have injective monodromy. By construction the fibers of the projection $W^{2}_b\to Z'_b$ are double coverings of $Z''_b$ ramified in the two points of $\left(\left\{z'\right\}\times Z''_b\right)\cap D^{1}$ for $z'\in Z'_b$ and thus have genus $2+8(g-1)$. Hence, the same holds for the fibers of the projection $W^{2}\to Z'$.

For a fixed $b_0\in B'$ the same arguments as in Section \ref{candfi} now show that we can apply Proposition \ref{propositioninjg} to the short exact sequence
\[
 \xymatrix{ 1\ar[r]& \pi_1(W^{2}_{b_0})\ar[r]& \pi_1(W^{2})\ar[r]& \pi_1(B')\ar[r]& 1,}
\]
with $f_1: W^{2}\to Z'$ and $f_2: W^{2}\to Z''$ the natural projections and $u_1:W^{2}_{b_0}\to Z'_{b_0}$ (resp. $u_2: W^{2}_{b_0}\to Z''_{b_0}$) the two distinct holomorphic submersions of the Kodaira fibration $W^{2}_{b_0}$. In particular, it follows that the monodromy of each of the fibrations $f_i$ is injective. 

To conclude the proof of Theorem~\ref{basearbitraire} we simply have to explain the irreducibility statement. We claim that the monodromy associated to the fibration $f_i$ ($i=1, 2$) is irreducible when restricted to any finite index subgroup of the fundamental group of the fiber of the fibration $Z'\to B'$ or $Z''\to B'$ (hence is also irreducible on any finite index subgroup of $\pi_{1}(Z')$). But this follows from Shiga's result saying that the monodromy of any holomorphic family of Riemann surfaces over a Riemann surface of finite type is irreducible, see~\cite{shiga}. We also refer the reader to~\cite{nicolasc} for another proof of this fact. 

\begin{rem} Of course, the irreducibility statement that we have just proved also holds in the context of Theorem~\ref{allve}. 
\end{rem}

\vspace{.5cm}

We now explain how to deduce the proof of Theorem~\ref{veofmcg} from the above construction. For $g\ge 2$, we denote by $\mathcal{T}_{g}$ (resp. $\mathcal{T}_{g,\ast}$) the Teichm\"uller space of closed Riemann surfaces of genus $g$ (resp. of genus $g$ with one marked point). We refer to~\cite{FarbMargalit} for basic facts on Teichm\"uller spaces. The spaces $\mathcal{T}_{g}$ and $\mathcal{T}_{g,\ast}$ are complex manifolds of complex dimension $3g-3$ and $3g-2$ respectively. The group $\mcg(S_g)$ (resp. $\mcg (S_g,\ast)$) acts properly discontinuously by holomorphic maps on $\mathcal{T}_{g}$ (resp. $\mathcal{T}_{g,\ast}$). In particular, the quotient of either $\mathcal{T}_{g}$ or $\mathcal{T}_{g,\ast}$ by any torsion-free subgroup of the corresponding mapping class group is also a complex manifold, while the moduli spaces $\mathcal{M}_{g}:=\mathcal{T}_{g}/\mcg(S_g)$ and $\mathcal{M}_{g,\ast}:=\mathcal{T}_{g,\ast}/\mcg (S_g,\ast)$ are complex orbifolds. It is well-known that $\mcg(S_g)$ and $\mcg (S_g,\ast)$ have torsion-free subgroups of finite index. A sequence of such torsion-free finite index subgroups is given by the level $m$ congruence subgroups $\mcg(S_g)\left[m\right]<\mcg(S_g)$ and $\mcg (S_g,\ast)\left[m\right]<\mcg (S_g,\ast)$ with $m\geq 3$, which are defined as the kernels of the natural morphisms to ${\rm Sp}(2g,\Z/m\Z)$ induced by the action of each mapping class group on the homology of the underlying (unmarked) surface. 
The \emph{universal family of curves of genus $g$} is defined by the holomorphic map $\mathcal{M}_{g,\ast}\to \mathcal{M}_{g}$ of complex orbifolds which ``forgets'' the marked point. It induces the Birman exact sequence 
{\small
\begin{equation}\label{birmanes}
 \xymatrix{ 1 \ar[r]& \pi_1(S_g) \ar[r]& \mcg(S_g,\ast)\ar[r]& \mcg(S_g)\ar[r]& 1}
\end{equation}}
\noindent on orbifold fundamental groups. Its monodromy is of course injective. Passing to the finite index subgroup $\mcg(S_g)[3]< \mcg(S_g)$ and its preimage $\mcg(S_g,\ast)[3]< \mcg(S_g,\ast)$, we obtain a short exact sequence
{\small
\begin{equation}
 \xymatrix{ 1 \ar[r]& \pi_1(S_g) \ar[r]& \mcg(S_g,\ast)[3]\ar[r]& \mcg(S_g)[3]\ar[r]& 1}
\end{equation}}
of torsion-free groups. It is induced by the surjective holomorphic map
$$ q: \mathcal{T}_{g,\ast}/\mcg(S_g,\ast)[3] \to \mathcal{T}_{g}/ \mcg(S_g)[3]$$ 
of complex {\it manifolds}. Its fibers are closed Riemann surfaces of genus $g$ and its monodromy is injective. We are now in the situation of Theorem \ref{basearbitraire} with $Z:=\mathcal{T}_{g,\ast}/\mcg(S_g,\ast)[3]$, $B:=\mathcal{T}_{g}/ \mcg(S_g)[3]$ and $\pi:=q$. Since by definition the group $\pi_1(Z) = \mcg(S_g,\ast)[3]$ is a finite index subgroup of $\mcg(S_g,\ast)$, this completes the proof of Theorem \ref{veofmcg}.

\begin{rem} There are probably plenty of morphisms as in Theorem~\ref{veofmcg}. Indeed, there are plenty of choices involved in our proof. Also, we have decided to perform our construction by taking {\it double} ramified covers, but one could also take cyclic ramified covers of higher orders, as in~\cite{Kod-67} or possibly use the constructions from~\cite{caupol}. It would be interesting to investigate how different all the morphisms $\Gamma \to \mcg (S_k)$ obtained in this way are (where $\Gamma$ is a finite index subgroup of $\mcg (S_g,\ast)$ and $k$ an integer).
\end{rem}

\begin{rem}\label{appholoeq} The morphisms that we construct naturally come with an equivariant holomorphic map $\mathcal{T}_{g,\ast}\to \mathcal{T}_{2+8(g-1)}$. This is interesting in view of the discusion in~\cite[\S 5]{arso-rigidity}. Conversely, let $\Gamma < \mcg (S_g,\ast)$ be a subgroup of finite index and let $f : \mathcal{T}_{g,\ast}\to \mathcal{T}_{k}$ be a holomorphic map, equivariant with respect to a representation of $\Gamma$ into $\mcg (S_k)$. Let $E_k\to \mathcal{T}_{k}$ be the universal family of curves of genus $k$ and let $f^{\ast}E_{k}\to \mathcal{T}_{g,\ast}$ be its pullback under $f$. By taking the quotient of $f^{\ast}E_{k}$ by $\Gamma \cap \pi_{1}(S_{g})$ (where $\pi_{1}(S_{g})$ is the normal subgroup appearing in Birman's short exact sequence~\eqref{birmanes}) we obtain a family of Kodaira fibrations over $\mathcal{T}_{g}$. Hence holomorphic equivariant maps $ \mathcal{T}_{g,\ast}\to \mathcal{T}_{k}$ are related to families of Kodaira fibrations with large deformation spaces.    
\end{rem}


\section{Further remarks}\label{furtherremarks}

This section contains two observations about residual properties and coherence of fundamental groups of surface bundles over surfaces. 
 
 \vspace{.4cm}
 
In~\cite{bregman} Bregman proves that certain Kodaira fibrations have a fundamental group which is not residually torsion-free nilpotent and asks whether there exist such fibrations with residually nilpotent (or residually torsion-free nilpotent) fundamental group. Here we observe that there are plenty of Kodaira fibrations with residually nilpotent fundamental group. This follows from Johnson's result discussed earlier, together with a result by Paris~\cite{paris}. This is also a result about surface bundles in general, which is independent of the consideration of complex structures. We summarize this observation in the following proposition. 
 
\begin{prop} Let $\pi : X\to B$ be any surface bundle over a surface. We assume that the monodromy of $\pi$ is injective. Then the fundamental group of $X$ is virtually residually $p$ and hence virtually residually nilpotent. 
\label{prop:virtnilpinjective}
\end{prop}

Note that any finite covering space of a surface bundle is again a surface bundle (with a base of higher genus), hence this provides plenty of examples of surface bundles with residually nilpotent fundamental group.

\vspace{.4cm}

\noindent {\it Proof.} Let $E\to X$ be the fibration whose fiber above a point $x\in X$ is the punctured surface $\pi^{-1}(x)-\{x\}$. By fixing a base point $x_0\in X$ and letting $F_{0}=\pi^{-1}(x_0)$ we obtain a morphism 
$$\varphi : \pi_{1}(X,x_{0})\to \mcg (F_{0},x_{0})$$
where $\mcg(F_{0},x_0)$ is the mapping class group of the once-punctured surface $(F_{0},x_{0})$. This group fits into Birman's short exact sequence:
$$\xymatrix{0 \ar[r] & \pi_{1}(F_0) \ar[r] & \mcg(F_{0},x_{0}) \ar[r] & \mcg(F_{0}) \ar[r] & 0}.$$
Using this short exact sequence, one sees easily that $\varphi$ is injective (this observation is classical; it relies on the fact that $\pi$ has injective monodromy). Since $\mcg(F_0,x_0)$ is virtually residually $p$~\cite{paris}, this completes the proof.\hfill $\Box$ 

As a consequence of Proposition \ref{prop:virtnilpinjective}, we also obtain the existence of surface bundles over surfaces with non-injective monodromy and virtually residually nilpotent fundamental group, as shown by the next remark.

\begin{rem}
 Let $\pi:X\to B$ be a surface bundle over a surface with injective monodromy and let $q:C\to B$ be a ramified covering with $q_{\ast}(\pi_1(C))=\pi_1(B)$ (e.g. a double covering with two ramification points of order 2).  
 Consider the base change  $Y=\left\{(x,c)\in X\times C \mid \pi(x)=q(c)\right\}$ of $X$ with respect to $q$. The projection $\pi_C: Y\to C$ equips $Y$ with the structure of a surface bundle over a surface and the map $Y\to X\times C$ induces an embedding $\pi_1(Y)\hookrightarrow \pi_1(X)\times \pi_1(C)$. In particular, it follows from Proposition \ref{prop:virtnilpinjective} that $\pi_1(Y)$ is virtually residually nilpotent. On the other hand, the monodromy of $\pi_C:Y\to C$ is not injective, since it is the composition of the epimorphism $q_{\ast}: \pi_1(C)\to \pi_1(B)$, which has non-trivial kernel, with the monodromy of $\pi:X\to B$.
\end{rem}

Finally, we close this section with an elementary observation about coherence of surface-by-surface groups, which is certainly well-known to experts (see e.g.~\cite{friedlvidussi2}). Recall that a group $G$ is called {\it coherent} if any finitely generated subgroup of $G$ is finitely presented.

\begin{prop} Let $\pi : X\to B$ be any surface bundle over a surface. If the corresponding monodromy morphism is not injective, then $\pi_{1}(X)$ contains a copy of a direct product of two nonabelian free groups and consequently is not coherent. 
\end{prop}
{\it Proof.} Let $R\lhd \pi_{1}(X)$ be the fundamental group of the fiber of $\pi$. By Lemma \ref{kmetpd} the centralizer $\Lambda$ of $R$ in $\pi_{1}(X)$ is a normal subgroup of $\pi_1(X)$ which is isomorphic to the kernel of the monodromy morphism
$$\pi_{1}(B)\to \out(R),$$
 and the subgroup generated by $R$ and $\Lambda$ is isomorphic to $R\times \Lambda$. If the monodromy is not injective, $\Lambda$ must be isomorphic to an infinite rank free group or a surface group. Hence, $R\times \Lambda$ (and also $\pi_{1}(X)$) contains a copy of $F_{2}\times F_{2}$. Non-coherence now follows from~\cite{baumslagroseblade}. 
 \hfill $\Box$ 

Let us explain why this observation seems interesting to us in relation to the study of the coherence of fundamental groups of aspherical K\"ahler surfaces. The study of coherence in the context of fundamental groups of complex surfaces was started by Kapovich in~\cite{kapovich98}. It was later pursued in~\cite{friedlvidussi,friedlvidussi2,kapovichgt,py16}. The outcome is that if $X$ is an aspherical K\"ahler surface with $b_{1}(X)>0$, then in most cases $\pi_{1}(X)$ is not coherent. More precisely, Friedl and Vidussi~\cite{friedlvidussi} prove that $\pi_{1}(X)$ is not coherent unless it is virtually the direct product of $\mathbb{Z}^{2}$ with the fundamental group of a closed surface or it is the fundamental group of a Kodaira fibration of virtual Albanese dimension $1$. In the last case, coherence is an open question. The proposition above shows for instance that the obvious examples of Kodaira fibrations with Albanese dimension $1$ (namely those whose monodromy has finite index in the mapping class group) have non-coherent fundamental group.

\bigskip
\bigskip
\begin{small}
\begin{tabular}{llll}
Claudio Llosa Isenrich & & & Pierre Py\\
Faculty of Mathematics & & & IRMA\\
Karlsruhe Institute of Technology & & & Universit\'e de Strasbourg \& CNRS\\
76131 Karlsruhe, Germany  & & & 67084 Strasbourg, France\\
claudio.llosa@kit.edu & & & ppy@math.unistra.fr\\    
\end{tabular}
\end{small}

\end{document}